\documentclass[12pt]{amsart}

\usepackage{amsmath,amsfonts,amsthm,amsopn}
\usepackage{graphicx}
\usepackage{epsfig,verbatim}
\usepackage{xspace}
\usepackage{mathrsfs}

\usepackage{xspace}
\newcommand{\supp}{\mathrm{supp}}
\newcommand{\spann}{\mathrm{span}}
\newcommand{\esssup}{\mathrm{esssup}}

\newcommand{\liL}{\lambda\in\Lambda}
\newcommand{\pil}{\pi (\lambda )}

\newcommand{\Lt}{\mathbf{L}^2 }
\newcommand{\Lz}{\mathbf{L}^2 (\mathbb{R}^d )}

\newcommand{\Lmp}{\mathbf{L}^{p,q}_m (\mathbb{R}^{2d} )}
\newcommand{\MS}{\mathbf{M}^{p,q}_m (\mathbb{R}^{d})} 
\newcommand{\Minf}{\mathbf{M}^{\infty} (\mathbb{R}^{d})} 
\newcommand{\Min}{\mathbf{M}^{\infty} } 

\newcommand{\SWD}{\mathcal{S}'(\mathbb{R}^d )}

\def\rd{\bR^d}

\def\rdd{{\bR^{2d}}}
\def\lrd{\mathbf{L}^2 (\mathbb{R}^d )}

\newcommand{\field}[1]{\mathbb{#1}}
\newcommand{\bR}{\field{R}} 
\newcommand{\bN}{\field{N}} 
\newcommand{\bZ}{\field{Z}} 
\newcommand{\bC}{\field{C}} 
\newcommand{\bQ}{\field{Q}} 

\newcommand{\tfa}{time-frequency analysis}

\newcommand{\stft}{short-time Fourier transform }

\newcommand{\tf}{time-frequency}

\newcommand{\fif}{if and only if}
\newcommand{\tfs}{time-frequency shift}

\newcommand{\modsp}{modulation space}

\def\cS{\mathcal{S}}

\def\cH{\mathcal{H}}

\def\cG{\mathcal{G}}

\def\cW{\mathcal{W}}

\def\cV{\mathcal{V}}

\def\vf{\varphi}
\def\inv{^{-1}}

\newtheorem{Th}{Theorem}

\newtheorem{Cor}[Th]{Corollary}

\newtheorem{Pro}[Th]{Proposition}
\newtheorem{Lem}[Th]{Lemma}
\theoremstyle{remark}
\newtheorem{rem}{Remark}

\newcommand{\ab}{\hspace{2pt}}
\newcommand{\g}{{\boldsymbol\varphi}}
\def\Wsp{{\boldsymbol W}} 

 \def\cS{\mathcal{S}}
 
 \def\cH{\mathcal{H}}

 \def\cG{\mathcal{G}}

 \def\cW{\mathcal{W}}

\def\rd{\bR^d}

\def\zd{\bZ^d}

\def\rdd{{\bR^{2d}}}
\def\zdd{{\bZ^{2d}}}

\def\lrd{\mathbf{L}^2(\rd)}

\def\Lpqm{\mathbf{L} ^{p,q}_m}
\def\Lmpq{\mathbf{L} ^{p,q}_m}
\def\lpqm{\ell_{m}^{p,q}}

\def\Mmpq{\mathbf{M}_m^{p,q}}
\newcommand{\bm}{\mathbf{M}}

\def\intrdd{\int_{\rdd}}

\newcommand{\mwgf}{multi-window Gabor frame}

\begin{document}
\title[Time-Frequency Partitions]{Time-Frequency Partitions and Characterizations of  Modulation
  Spaces with Localization Opertors}
\author{Monika D\"orfler and Karlheinz Gr\"ochenig}
\address{Institut f\"ur Mathematik, Universit\"at Wien, Alserbachstrasse 23
  A-1090 Wien,   Austria }
\email{\{monika.doerfler,karlheinz.groechenig\}@univie.ac.at} 
\thanks{M.\ D.\ was supported by the FWF grant T~384-N13,  K.\ G.\ was
  supported by the  Marie-Curie Excellence Grant 
  MEXT-CT-2004-517154 and in part by the  National Research Network S106 SISE of the
Austrian Science Foundation (FWF)}

\begin{abstract}
We study families of  \tf\ localization operators and derive a
new characterization  of modulation spaces. This characterization
relates the size of the localization operators  to
the global \tf\ distribution. 
As a by-product,  we obtain a new proof for the existence of  multi-window
Gabor frames and extend the structure theory of Gabor frames. 
\end{abstract}
\date{\today}
\maketitle

\section{Introduction}\label{int}

A \tf\ representation transforms  a function $f$ on $\rd $ into 
a function on the \tf\ space   $\rd \times \rd $. The goal is to
obtain a description of $f$ that is local both in time and in
frequency~\cite{brka03,gora98}.  The  standard \tf\ representations,
such as the \stft\ and its various modifications known as  Wigner 
distribution,  radar ambiguity function, Gabor transform,  all 
encode   \tf\ information.    However, the pointwise interpretation
of such a \tf\ representation meets difficulties because,  by the uncertainty
principle,   a small region in the \tf\ plane does not possess a
physical meaning. Therefore  the question arises in which sense the
short-time Fourier transform  describes the local properties of a function and its
Fourier transform. 

Following Daubechies~\cite{Daube88}, we use \tf\
localization operators to give meaning to the local \tf\ content. By
investigating a whole family of localization operators and glueing
together the local pieces, we are able to characterize the global
\tf\ distribution  of a function. In more technical terms, our main
result provides a new characterization of modulation spaces.

We define the \stft\   (STFT)  of
a function $f\in \mathbf{L}^2 (\mathbb{R}^d) $ with respect to a
window function $\varphi\in \mathbf{L}^2(\mathbb{R}^d) $  as
\begin{equation}
  \label{Def:STFT}
  \mathcal{V}_{\varphi} f(x,\omega) = \int _{ \bR^d } f(t) \bar{\varphi} (t-x) e^{-2\pi
  i \omega \cdot t} \, dt 
,\mbox{ for all } z = (x,\omega) \in \mathbb{R}^{2d}.
\end{equation} 

  The STFT $\mathcal{V}_{\varphi} f(z)$
  is a  measure of  the \tf\ content near the point $z$
in the \tf\ plane $\rdd $. However, the  STFT 
cannot be  supported on a set  of finite measure by results in
~\cite{jaming98,janssen98,wilczock98}. This fact  complicates the interpretation of local
information obtained from the STFT. In particular, it is impossible to
construct a 
projection operator that satisfies $V_\varphi
(P_\Omega f) = \chi _\Omega \cdot V_\varphi f $.  As a remedy   one resorts to the following definition
of localization operators.

  
We  denote  translation operators by $T_xf(t) = f(t-x)$ and  time-frequency shifts by $\pi(z)f(t)  = e^{2\pi i \omega
  \cdot t} f(t-x)$ for $x,\omega,t\in \rd $.
Fix a non-zero function 
  $\varphi\in\mathbf{L}^2(\mathbb{R}^{2d})$ (a so-called  window
  function)   and a symbol $\sigma \in \mathbf{L}^1(\rdd )$. Then
  the {\it time-frequency localization operator}
$H_{\sigma}$ acting on a function $f$  is 
defined as
\[H_{\sigma} f = \int_{\mathbb{R}^{2d}} \sigma (z) \mathcal{V}_{\varphi} f(z
)\pi (z) \varphi\ab dz =
\mathcal{V}_{\varphi}^{\ast}\sigma\mathcal{V}_{\varphi}  f.\]
The integral is defined strongly on   many function spaces, in
particular on $\mathbf{L}^2(\rd )$. 
  A useful alternative  definition of $H_\sigma$ is the weak definition 
\begin{equation}
  \label{eq:f1}
  \langle H_\sigma f,g\rangle _{L^2(\rd )}  = \langle \sigma
  \mathcal{V}_{\varphi} f, \mathcal{V}_{\varphi} g\rangle _{L^2(\rdd
    )} \, .
\end{equation}
This definition can be easily extended to distributional symbols
$\sigma \in \cS ' (\rdd ) $. 
The  subtleties of the
definition and boundedness properties between various spaces have been
investigated in many papers, see \cite{CG03,Toft04,wo02} for a sample of
results. 
%

If $\sigma$ is non-negative and has compact
support  in $\Omega \subseteq \rd$, then  $H_{\sigma}f$ can be
interpreted as the part of $f$ that  lives essentially on    $\Omega$
in the time-frequency plane, and so  $H_\sigma $ may be taken as  a
substitute for the non-existing projection  onto the region
$\Omega $ in the \tf\ plane. 
 
In this paper 
 we investigate the behavior of  an entire collection of localization
 operators. Namely, given a lattice $\Lambda \subseteq \rdd $ of the
 \tf\ plane, we consider the collection of operators 
$\{ H_{T_\lambda \sigma} :\lambda\in \Lambda\}$ and the   mapping $f \to
\{H_{T_\lambda\sigma }f\}$. If the supports of $T_\lambda \sigma $
cover $\rdd $, then $\{H_{T_\lambda\sigma }f, \lambda \in \Lambda \}$
should contain enough information to recover $f$ from its local
components. In particular, the set $\{H_{T_\lambda\sigma }f: \lambda
\in \Lambda\}$ should carry the complete information about the global
\tf\ properties of $f$.  We  make this intuition precise  and
derive  a new characterization of {\it \modsp s} from it. Similar to Besov spaces,
modulation spaces are smoothness spaces, but the smoothness is
measured by means of time-frequency distribution rather than
differences and derivatives.  Here, we establish a
correspondence between the behavior of the sequence $\|H_{T_\lambda
  \sigma }f\|_2, \lambda \in \Lambda $, and the membership of $f$ in a
\modsp . 

As a special case of our main theorem we formulate the following
result. 

\begin{Th}\label{introd}
Fix a non-zero function    $\vf $ in the Schwartz space $ \cS
(\rd ) $ and a  weight function $m$ on $\rdd $
that satisfies $m(z_1+z_2) \leq C (1+|z_1|)^N m(z_2)$ for some
constants $C,N\geq 0$ and all $z_1,z_2\in \rdd $. Then a tempered
distribution $f$ satisfies 
\begin{equation}
  \label{eq:c11}
  \Big( \intrdd  |\cV _\vf f (z )|^p m(z)^p \, dz \Big)^{1/p} <\infty
  \, , 
\end{equation}
\fif\ 
\begin{equation}
  \label{eq:c12}
  \Big(\sum_{\lambda\in\Lambda}
  \|H_{T_{\lambda} \sigma} f\|_2^p m(\lambda )^p \Big)^{1/p} <\infty \, .
\end{equation}
 \end{Th}
The expression in~\eqref{eq:c11} is just the norm of $f$ in the
\modsp\ $\mathbf{M}^p_m(\rd )$. Our main result shows that the  expression
in~\eqref{eq:c12} (using the \tf\ components of $f$)  is an equivalent norm on the
\modsp\ $\mathbf{M}^p_m(\rd )$. 

In pseudodifferential calculus one often defines spaces by conditions
on their time-frequency components. For instance, Bony, Chemin, and
Lerner~\cite{boch94,BoL} introduced a Sobolev-type space $H(m)$  by using Weyl
operators instead of localization operators. For the (extremely
simplified) case of a constant Euclidean  metric on the \tf\ plane, a
distribution $f$ belongs to $H(m)$, whenever for some test function
$\psi$ on $\rdd $
\begin{equation}
  \label{eq:f3}
\|f\|_{H(m)}^2 = \intrdd \| (T_Y\psi )^w f\|_2^2 \, m(Y) \, dY \, ,  
\end{equation}
is finite, 
where  $\sigma ^w$ is the
Weyl operator corresponding to the symbol $\sigma $.  The only difference
between \eqref{eq:f3} and  \eqref{eq:c12} is
the use of Weyl calculus instead of \tf\ localization operators and a
continuous definition instead of a discrete one. It was understood
only recently that $H(m)$ coincides with the \modsp\ $M^2_m(\rd )$ and
that  \eqref{eq:f3} is an equivalent norm on $M^2_m(\rd
)$~\cite{GT09}. Thus Theorem~\ref{introd} can be interpreted as an
extension of \cite{boch94} to $L^p$-like spaces. 
 
Let us also mention that in the language of ~\cite{su06}, the
operators   $\{H_{T_{\lambda} \sigma}
,\,\lambda\in\Lambda\}$ form a g-frame for  $\mathbf{L}^2 (\bR^d
)$. Our construction seems to be one of the few non-trivial examples
of $g$-frames that are not frames.

In this paper we prove the norm equivalence of Theorem~\ref{introd} for a large
class of \modsp s and arbitrary \tf\  lattices. For a  rather restricted class of lattices, namely lattices with integer oversampling,   an analogous result
was derived in ~\cite{dofegr06} for  unweighted \modsp
s. The main arguments for the integer lattice were based on Zak
transform methods and interpolation. For a general lattice, these
methods are no longer 
available, and we have to develop a completely new  approach to some of the
key arguments. 

As a by-product of the new techniques we have found several 
results of independent interest. 

\begin{itemize}
\item We formulate several  structural results and characterizations
  of Gabor frames for  multi-window Gabor frames. 
\item We prove  a finite intersection property for time-frequency
  invariant subspaces of the distribution space
  $\mathbf{M}^\infty(\mathbb{R}^d)$. This property resembles the
  finite intersection property that characterizes compact sets. 
\item We give a new, independent proof for the  existence of
  multi-window Gabor frames with well-localized windows. Previous
  proofs were based on coorbit theory~\cite{fg92chui} and the theory of
  projective modules~\cite{lu09}. Our proof provides additional insight
  how the windows can be chosen.  
\item We derive precise  estimates for the localization  of the
  eigenfunctions of a localization operator.
\end{itemize}

This paper is organized as follows. In  Section~2  we recall
necessary  facts from time-frequency analysis. On the one had,  we
introduce modulation spaces and explain their characterization by
means of multi-window Gabor frames. 
On the other hand, we 
state and
prove several  properties of localization operators. In
Section~\ref{Se:Main}, we formulate and prove our main result
(Theorem~\ref{Th:Main}). In Section~\ref{Se:Con} we analyze some of the 
consequences of Theorem~\ref{Th:Main} and its proof. In the appendix
we collect and sketch the proofs of some of the structural results on
Gabor frames. 

\section{Time-Frequency Analysis of Functions and Operators}
\subsection{Modulation Spaces}\label{Se:ModSp} 
Modulation spaces are a class of  function spaces associated to the
short-time Fourier transform~\eqref{Def:STFT}. Note that  for  a  suitable test function $ \varphi $, the
\stft\ can be extended to distribution spaces by duality and
$\mathcal{V}_\varphi f (z) = \langle f, \pi (z) \varphi \rangle $.

For the  standard definition of modulation spaces,   we  fix a non-zero "window
function" $g \in \cS (\rd )$   and  consider   moderate   weight
functions $m $ of polynomial growth, i.e.,  $m$ satisfies  $m (z_1+z_2) \le C  (1+|z_1|)^s m
(z_2)$, $z_1,z_2 \in \rdd $ for some  $C,s \ge 0$. 
Given a moderate weight
$m$ and $1\leq p,q \leq \infty $,  the modulation
space $\MS$ is defined as the space of all tempered distributions
$f\in\SWD$ with $\mathcal{V}_g f\in \Lmp$, with  norm
\begin{equation}
  \label{eq:moddefc}
\|  f \|_{\MS} = \|\mathcal{V}_g f\|_{\Lmp}\, .
\end{equation}
If $p=q$, we write $\mathbf{M}^p_m(\mathbb{R}^d)$.\\

For weight functions of faster growth we have to resort to different
spaces of test functions and distributions. Let $g(t) = e^{-\pi t\cdot
  t}$ be the Gaussian window and $\cH _0 = \mathrm{span}\, \{ \pi (z)
g : z\in \rdd \}$ be the linear space of all finite linear
combinations of \tfs s of the Gaussian. Let $\nu $ be a
submultiplicative even weight function on $\rdd $  and $m$ be a $\nu $-moderate
function; this means that $\nu (z_1+z_2) \leq \nu (z_1) \nu (z_2)$,
$\nu (z) = \nu (-z)$ and
$m(z_1+z_2) \leq \nu (z_1) m (z_2)$ for all $z, z_1,z_2\in \rdd $. For
$1\leq p,q<\infty $ the \modsp\ $\Mmpq (\rd )$ is then defined as the
closure of $\cH _0$ in the norm $\|f\|_{\MS }$ as in
\eqref{eq:moddefc}. If $p=\infty $ or $q=\infty $, we take a
weak-$^*$-closure of $\cH _0$. These general \modsp s possess the
following properties. Assume that  $m$ is $\nu $-moderate  and $1\leq p,q\leq \infty $, then 
\begin{equation}
  \label{eq:c19}
  \bm _\nu ^1 (\rd )  \subseteq \Mmpq (\rd ) \subseteq M^{\infty
  }_{1/\nu }(\rd ) =   \bm _\nu ^1 (\rd ) ^* \, .
\end{equation}

Further,  if $\vf \in   \bm _\nu ^1 (\rd )$, then 
\begin{equation}
  \label{eq:c20}
  \|\cV _\vf f \|_{\Lpqm } \asymp \|\cV _g f\|_{\Lpqm } = \|f\|_{\Mmpq
  } \, .
\end{equation}
thus different windows in $  \bm _\nu ^1 (\rd )$ yield equivalent
norms on $\Mmpq $. 

The embedding \eqref{eq:c19} says that  $  \bm _\nu ^1 (\rd )$ may
serve as a space of test functions and $\bm ^\infty _{1/\nu } (\rd )$
as a space of distributions for all \modsp s $\Mmpq $ with a $\nu
$-moderate weight $m$. 

If $\nu _s (z) = (1+|z|)^s, s\geq 0$ and $m$ is $\nu _s$-moderate, then we have 
$$
\cS (\rd ) \subseteq   \bm _{\nu _s} ^1 (\rd ) \subseteq \Mmpq (\rd )
\subseteq \bm ^\infty _{1/\nu _s }(\rd ) \subseteq \cS ' (\rd ) \, ,
$$
in agreement with the standard definition, but for $\nu (z) = e^{a
  |z|^b}$ with $a>0$ and $0<b\leq 1$ we have 
$$
  \bm _\nu ^1 (\rd ) \subseteq \cS (\rd ) \subseteq \cS ' (\rd )
  \subseteq \bm ^\infty _{1/\nu }(\rd ) \, .
$$
In the sequel we will start with a submultplicative weight $\nu $ and
take $\bm ^\infty _{1/\nu }(\rd )$ as the appropriate distribution
space. Our results hold for arbitrary  submultiplicative weights
$\nu $.

For the detailed theory of  \modsp s we refer to \cite[Ch.
11--13]{gr01}, for a discussion of weights and possible distribution
spaces see \cite{gro07c}.

\textbf{Sequence space norms.} Recall that  a time-frequency lattice $\Lambda$ is a discrete subgroup of $\bR^{2d}$
of the form $\Lambda = A\bZ^{2d}$ for some invertible real-valued
$2d\times 2d$-matrix $A$.\\
   Given a lattice $\Lambda \subseteq \rdd
$ with relatively compact fundamental domain $Q$, the   discrete  space $\ell^{p,q}_m
(\Lambda )$ consists  of all sequences $\mathbf{a} =
(a_{\lambda })_{\lambda \in \Lambda}$ for which the norm
\begin{equation}
  \label{eq:c21}
  \|\mathbf{a} \|_{\lpqm } = \|\sum _{\lambda \in \Lambda } |a_\lambda
  | \chi _{\lambda +Q} \|_{\Lmpq } 
\end{equation}
  is finite. If $\Lambda = a\zd \times b\zd $, then this definition
  reduces to the usual mixed-norm space $\lpqm (\zdd )$ with norm 
\[\|\mathbf{a}\|_{\ell^{p,q}_m }
=\Big(\sum_{n\in\mathbb{Z}^d}\Big(\sum_{k\in\mathbb{Z}^d}|a_{kn}|^p
m(ak,bn)^p\Big)^{q/p}\Big)^{1/q}\, .\]

As a technical tool we will need  amalgam
spaces (in one place only).  A
measurable function $F$ on $\mathbb{R}^{2d}$ belongs to the (Wiener)
amalgam space $\Wsp (\mathbf{L}^{p,q}_m )$, if the sequence of local
suprema 
\[a_{kn} = \esssup_{x,w\in [0,1]^d} |F(x+k, \omega +n ) | = \|F\cdot T_{(k,n)}\chi\|_\infty\]
belongs to $\ell^{p,q}_m (\mathbb{Z}^{2d})$. The norm on $\Wsp
(\mathbf{L}^{p,q}_m )$ is $\| F\|_{\Wsp (\mathbf{L}^{p,q}_m )} = \|
a\|_{\ell^{p,q}_m }$.  See~\cite{he03} for an introductory
article.
We need  their
 behavior under convolution and their properties under sampling. 
\begin{itemize}
\item[(a)] Convolution in Wiener amalgam spaces: Let $1\leq p,q\leq\infty$ and let $m$ be a $\nu$-moderate weight. Then 
\begin{equation}\label{Eq:WAsamp}
\|F\ast G\|_{\Wsp (\mathbf{L}^{p,q}_m )}\leq C \|F\|_{\Wsp (\mathbf{L}^{p,q}_m )}\|G\|_{\mathbf{L}^1_\nu}.
\end{equation}
\item[(b)] Sampling in Wiener amalgam spaces:
For $F\in \Wsp (\mathbf{L}^{p,q}_m)$ the following sampling property holds: 
\begin{equation}\label{Eq:lpsamp0}
\|F|_\Lambda\|_{\ell^{p,q}_{{m}}} \leq C_{\Lambda}\| F\|_{\Wsp (\mathbf{L}^{p,q}_m)}.
\end{equation} 
\end{itemize}
These statements are proved in~\cite{he03} or~\cite
[Prop.~11.1.4., Thm.~11.1.5.]{gr01}.
\subsection{Gabor frames}\label{Se:MWGF} 
Gabor
frames are closely  linked to modulation spaces. They constitute
``basis-like'' sets for \modsp s and are used to 
 characterize the membership in a \modsp\  by the magnitude  of
coefficients in the corresponding series expansion. 

For a given lattice $\Lambda\subseteq \rdd $ and a  window function
$\varphi\in\mathbf{L}^2 (\mathbb{R}^d )$, let $\mathcal{G} (\varphi,
\Lambda )$ denote  the  set of functions 
$\{ \pi (\lambda )\varphi : \lambda \in \Lambda\} $
 in $\mathbf{L}^2 (\mathbb{R}^d )$.
The operator
\[S_{\varphi }f  = \sum_{\lambda\in\Lambda}\langle f,
\pi (\lambda )\varphi\rangle  \pi (\lambda )\varphi \] is the {\it
frame operator} corresponding to  $\mathcal{G} (\varphi, \Lambda  )$.
If $S_\vf $ is bounded and invertible on $\lrd $, then $\mathcal{G} (\varphi, \Lambda
)$ is called  a
\emph{Gabor frame} for  $\mathbf{L}^2 (\mathbb{R}^d )$. 
This property is equivalent to the existence of two constants $A,B >0$
such that 
\begin{equation}
  \label{eq:c13}
  A \|f\|_2^2 \leq \sum _{\lambda \in \Lambda } |\langle f, \pi
  (\lambda )g\rangle |^2 = \langle S_\vf f, f\rangle \leq B\|f\|_2^2 \qquad \text{for all } \,
  f\in \lrd \, .
\end{equation}

Using  several windows $\g = (\varphi _1, \dots \vf _n)$, we say that  the union
$\bigcup_{j=1}^n\mathcal{G} (\varphi_j, \Lambda )$ is
a   \emph{multi-window Gabor frame},  if  the  associated 
frame operator  given by
\begin{equation}\label{Def:FOP}
 S_\g f  = \sum_{j=1}^n \sum_{\liL}\langle f,
\pi (\lambda )\varphi_j\rangle  \pi (\lambda )\varphi_j =
\sum_{j=1}^n S_{\varphi_j} f	
\end{equation}
is invertible on $\lrd $. 
 The frame operator can be expressed as  the composition  of the  analysis operator   
$C_{\boldsymbol\varphi,\Lambda}$  defined by 
     $$
C_{\boldsymbol\varphi,\Lambda}(f)(\lambda , j )  = \langle
     f,\pil \varphi_j\rangle  \, ,  \qquad \lambda \in \Lambda ,
     j=1,\dots ,n \, .
$$
 and the  synthesis operator   
  $D_{\boldsymbol\varphi,\Lambda}$ defined  by 
      $D_{\boldsymbol\varphi,\Lambda}(\mathbf{c} ) =
      \sum_{\liL}\sum_{j = 1}^n c_{\lambda ,j} \pi (\lambda )
      \varphi_j $. Then       $S_{\boldsymbol\varphi,\Lambda} =
      D_{\boldsymbol\varphi,\Lambda}\circ
      C_{\boldsymbol\varphi,\Lambda}$.

\subsection{Characterization of Modulation Spaces with  Gabor Frames}\label{Se:charMSMW}
The following  characterization of modulation spaces by means of multi-window
Gabor frames  is a  central result in \tfa\ and useful in many
applications. It is  crucial   for  the
proof of our main theorem (Theorem~\ref{Th:Main}).

\begin{Th}
 \label{fund}
Let  $\nu $  be a submultiplicative weight  on $\rdd $ satisfying the
condition $\lim _{n\to \infty } \nu (nz)^{1/n}= 1$ for all $z\in \rdd
$ and let  $m$ be a $\nu $-moderate weight and $1\leq p,q \leq \infty
$. 
 Assume further  that  $\bigcup_{j=1}^n\mathcal{G} (\varphi_j,
  \Lambda )$ is 
a   multi-window Gabor frame and that $\vf _j \in \mathbf{M}^1_\nu (\rd )$ for
$j=1, \dots , n$.
 
(i)   A distribution $f$ belongs to $\mathbf{M}^p_m(\rd )$, if and
only if $C_{\vf _j} f \in \ell^{p}_{{m}}$ for $j=1, \dots , n$. 
In this case there exist constants $A,B>0$,  such that, for all  $
f\in \mathbf{M}^p_m(\rd )$, 
\begin{equation*}\label{Eq:MuWiEquiv}
A\| f\|_{\mathbf{M}^p_m}\leq \Big( \sum_{\lambda\in\Lambda}\Big(\sum_{j=1}^n |\langle f,\pi (\lambda )
\varphi_j\rangle |^2\Big)^{p/2} m(\lambda )^p \Big)^{1/p} \leq B\|
f\|_{\mathbf{M}^p_m} \, .
\end{equation*}

(ii) Assume in addition that $\Lambda = a\zd \times b\zd $ is a
separable lattice. Then a   distribution $f$ belongs to $\MS $ if and
only if each sequence $C_{\varphi_j}f(ak,bl) = \langle f,\pi (ak,bl)
\varphi_j\rangle$ belongs to $\ell ^{p,q}_{m } (\zdd )$. 
In this case there exist constants $A$ and $B$ depending on $p,q,m$ such that, for all $ f\in \Mmpq$
\begin{equation}\label{Eq:MuWiEquivb}
A\| f\|_{\Mmpq }\leq
\sum_{l\in\mathbb{Z}}\Big(\sum_{k\in\mathbb{Z}}\Big(\sum_{j=1}^n |
\langle f,  \pi (ak,bl) \vf _j \rangle  |^2\Big)^{p/2}
m(ak,bl)^p\Big)^{q/p}\Big)^{1/q}  \leq B\| 
f\|_{\Mmpq }
\, .
\end{equation}

(iii) Let  $\Lambda \subseteq \rdd $ be an arbitrary lattice
and $Q$ be a relatively compact fundamental domain of $\Lambda $. Then
a distribution $f$ belongs to
$\MS $, \fif\ the function $\sum _{\lambda \in \Lambda } \Big(\sum _{j=1}^n |\langle f, \pi (\lambda )
\vf _j \rangle |^2\Big)^{1/2}  \chi _{\lambda +Q} $ belongs to $ \Lmpq (\rdd )$. In
this case there exist constants $A,B>0$, such that,  for all  $f\in \Mmpq (\rdd )$,
$$
A\|f\|_{\Mmpq} \leq \|\sum _{\lambda \in \Lambda } \Big(\sum _{j=1}^n |\langle f, \pi (\lambda )
\vf _j \rangle |^2\Big)^{1/2}  \chi _{\lambda +Q} \|_{\Lmpq }
\leq B \|f\|_{\Mmpq } \, .
$$
 \end{Th}

Note that (ii) follows from (iii), since  for $Q=[0,a]^d\times
[0,b]^d$  the norm equivalence $\|\sum _{k,l\in \zdd } a_{kl}
\chi _{(ak,bl)+Q} \|_{\Lmpq } \asymp \|a\|_{\lpqm } $ holds . 


Theorem~\ref{fund} has a long history. It  extends  the basic characterizations of
\modsp s by Gabor frames to \mwgf s.   For Gabor frames with a single
window and  lattices of the form  $\Lambda = a\zd \times b\zd $
with $ab \in \bQ $ Theorem~\ref{fund} was proved in \cite{fg97jfa}. For general
lattices it follows from the main  result in  \cite{GL03} and the
techniques in \cite{fg97jfa}. See also the discussion 
  in~\cite[Ch.~13]{gr01}. 
The proofs for  multi-window Gabor frames require only few modifications,
we therefore  postpone a discussion to the appendix. 


\subsection{A New Characterization of  Multi-Window Gabor Frames}\label{Se:charMWGF}
The proof of our main statement relies on a characterization of \mwgf
s without using inequalities. The following lemma is a generalization
of  ~\cite{gr07} from Gabor frames to  multi-window Gabor  frames. 
\begin{Lem}\label{Le:MuwiwithoutIneq}Assume that  $\varphi_j\in
  \mathbf{M}^1 (\rd )$ for  $j = 1,\ldots , n$. Then the following
  properties are equivalent. 
  
\begin{itemize}
 \item[(i)] $\bigcup_{j=1}^n\mathcal{G} (\varphi_j, \Lambda )$ is a  multi-window Gabor frame for $\Lz$.
     \item[(ii)] The analysis operator $C_{\boldsymbol\varphi,\Lambda}$ is one-to-one from $\mathbf{M}^{\infty} (\bR^d )$ to $\ell^{\infty } (\Lambda , \bC^n )$.
 \end{itemize}
\end{Lem}
The idea of the proof will be given in Appendix~\ref{App}, where we 
will also list many more equivalent conditions. 

\subsection{Properties of Localization Operators}\label{Se:PropLoc}
We next recall  some elementary properties of the localization
operators $H_{T_\lambda \sigma }$. 
Time-frequency localization operators have been introduced
and studied by Daubechies~\cite{daubechies90,Daube88} and 
 Ramanathan and Topiwala~\cite{RT93}, and are  also called STFT multipliers, \tf\ Toeplitz
operators, Wick operators, \tf\ filters, etc. 
They are a popular tool in 
 signal analysis for \tf\ filtering or nonstationary filtering~\cite{lo07-2,ro98}, in
quantization procedures in physics~\cite{Berezin71},  or in the approximation
of pseudodifferential operators~\cite{CF78,Lerner}. For a detailed
account of the early theory we refer to Wong's book~\cite{wo02}, for a  study of boundedness and Schatten
class properties to~\cite{CG03,cogrro06,fega06,Toft04}. 

\begin{Lem}[Intertwining property] \label{Le:TFshiftHl}
 If $\sigma \in \mathbf{L}^\infty(\mathbb{R}^{2d})$, $\varphi \in \Lt(\mathbb{R}^d) $, and $\lambda\in \Lambda $, then
\[\pi (\lambda)\ab H_{\sigma}\ab \pi (\lambda)^{\ast} = H_{T_\lambda \sigma}.\]
\end{Lem}
The proof consists of a simple calculation, see~\cite[Lemma~2.6]{dofegr06}.\\

For estimates of the STFT of $H_\sigma f$ we introduce the formal
adjoint of $\cV _\vf $, namely
$$
\mathcal{V}_{\varphi}^{\ast} F = \int_{\mathbb{R}^{2d}} F(z )\pi
(z) \varphi \ab dz \, ,
$$
which maps functions on $\rdd $ to functions or distributions on $\rd
$. With this notation we can write the localization operator
$H_\sigma$  as
$$
H_\sigma f = \cV ^*_{\vf } \big( \sigma \, \cV _\vf f \big) \, .
$$
The STFT of $\cV ^*_{\vf } F $ satisfies a fundamental pointwise estimate~\cite[Proposition 11.3.2.]{gr01}:
\begin{equation}\label{eq:STFTconva}
|\mathcal{V}_{\varphi} ( \mathcal{V}_{\varphi}
^{\ast}F ) (z)| \leq
\Big(|\mathcal{V}_{\varphi} \varphi|\ast|    F|\Big) (z) \quad \quad
\forall z\in \rdd \, .
\end{equation}
We note that for $F = \sigma \cV _\vf f$ this estimate becomes
\begin{equation}
  \label{eq:c14}
  |\cV _\vf (H_\sigma f)(z)| = | \cV _\vf \Big(\cV _\vf ^* (\sigma \cV
  _\vf f)\Big)(z)| \leq \Big( |\cV _\vf \vf | \ast (\sigma \, |\cV
  _\vf f|)\Big) (z) \, .
\end{equation}
Thus the \stft\ of $H_\sigma $ is a so-called product-convolution
operator. The standard boundedness results for localization operators
can be easily deduced from the well established results for product
convolution operators~~\cite{busby-smith81}.

 Estimate~\eqref{eq:c14} is quite  useful for the derivation of norm
estimates.
In the following, we fix a non-negative symbol $\sigma$ and investigate the set  of operators $\{H_{T_\lambda \sigma} :\lambda\in\Lambda\}$. To simplify notation we will write $H_\lambda$ instead of
$H_{T_\lambda\sigma}$, and sometimes  $H_0 = H_\sigma$ by some abuse of notation. 

\begin{Lem}\label{Le:PropH_0}
(i) Assume that  $\sigma \in \mathbf{L}^1 (\mathbb{R}^{2d})$, $\sigma
\geq 0 $  and that $\varphi \in \Lt  (\mathbb{R}^d)$. Then each
 $H_\lambda$, $\lambda \in\Lambda$,  is a positive trace-class operator.

(ii) If,  in addition, $\varphi \in \mathbf{M}^1_\nu (\rd )$ and
$\sigma \in
\mathbf{L}^1_\nu(\mathbb{R}^{2d})$, then each $H_\lambda$ is bounded 
from  $\mathbf{M}^\infty(\rd )$ into $\mathbf{M}^1_\nu (\rd )$.
 In particular, all eigenfunctions $\varphi _j $ of $H_\sigma $ belong to $\mathbf{M}^1_\nu (\rd )$.
 
 (iii)
Furthermore, if $\varphi \in \mathbf{M}^1_\nu (\rd )$ and
$\sigma \in
\mathbf{L}^1_\nu(\mathbb{R}^{2d})$, then each $H_\lambda$ is bounded 
from  $\mathbf{M}^\infty _{1/\nu }(\rd )$ into $\mathbf{L}^2(\rd )$. 
 \end{Lem}

\begin{proof}
Statement (i) is well-known, see, e.g., \cite{BC02,FN03,wo02}.

 To show (ii), we   use~\eqref{eq:c14}
to  obtain,  for $f\in \mathbf{M}^\infty (\rd )$, 
\begin{align}
  \|H_{\sigma} f\|_{\mathbf{M}^1_\nu}
  =&\,\, \|\cV _\vf (H_\sigma f)\|_{\mathbf{L}_\nu ^1} = \|\mathcal{V}_{\varphi}\mathcal{V}_{\varphi}
  ^{\ast}(\sigma\mathcal{V}_{\varphi} f )\|_{\mathbf{L}^1_\nu}\notag \\
 \leq& \,\,\big\| |\mathcal{V}_{\varphi} \varphi|\ast
 | \sigma \, \mathcal{V}_{\varphi} f |\big\|_{\mathbf{L}^1_\nu} \label{regineq} \\
 \leq&\,\, \|\mathcal{V}_{\varphi} \varphi\|_{\mathbf{L}^1_\nu} \, 
 \|\sigma\, \mathcal{V}_{\varphi} f \|_{\mathbf{L}^1_\nu} \, ,         \notag
\end{align}where we have used Young's inequality.
Since $\varphi \in \mathbf{M}^1_\nu (\rd )$ \fif\  $V_\varphi \varphi \in \mathbf{L}^1_\nu(\rdd )$
by \cite[Prop.~12.1.2]{gr01}, we find that
$$
  \|H_{\sigma} f\|_{\mathbf{M}^1_\nu} \leq \|\mathcal{V}_{\varphi}
  \varphi\|_{\mathbf{L}^1_\nu} \,  \|\sigma \|_{\mathbf{L}^1_\nu}
  \|V_{\vf } f\|_{\mathbf{L}^\infty }  \leq  C \|\sigma \|_{\mathbf{L}^1_\nu}
  \|f\|_{\mathbf{M}^\infty } \, , 
$$
and thus  $H_\sigma$ is bounded from $\mathbf{M}^\infty (\rd )$ to $\mathbf{M}^1_\nu (\rd )$.

The proof of (iii) is similar. Again, we apply~\eqref{eq:c14} to obtain
 for $f\in \mathbf{M}_{1/\nu}^\infty (\rd )$:
 $$ \|H_{\sigma} f\|_{\mathbf{L}^2}
  \leq
  \big\| |\mathcal{V}_{\varphi} \varphi|\ast
 | \sigma \, \mathcal{V}_{\varphi} f |\big\|_{\mathbf{L}^2} 
 \leq\, \|\mathcal{V}_{\varphi} \varphi\|_{\mathbf{L}^2} \, 
 \|\sigma\, \mathcal{V}_{\varphi} f \|_{\mathbf{L}^1} \, ,  
$$
 Hence, the result follows from
$$
\|\sigma\, \mathcal{V}_{\varphi} f \|_{\mathbf{L}^1}
= \int_{\mathbb{R}^{2d}}\sigma(z) |\mathcal{V}_{\varphi} f (z)|\nu(z)\frac{1}{\nu (z)}\, dz\\\notag
\leq  \|\sigma \|_{\mathbf{L}^1_\nu}\|f\|_{\mathbf{M}^\infty _{1/\nu
  }}\notag \, .
$$
\end{proof}
The spectral theorem for compact self adjoint operators provides  the following spectral representation of $H_\lambda$.

\begin{Cor}\label{Cor:Specrep}Assume 
$\varphi \in \mathbf{M}^1_\nu (\rd )$ and $\sigma \in \mathbf{L}^1_\nu(\mathbb{R}^{2d})$. Then there exists  a positive  sequence of
  eigenvalues $  \mathbf{c} = (c_j)\in \ell ^1$  and an
orthonormal system of eigenfunctions $\varphi_j\in \mathbf{M}^1_\nu (\rd )$, such that
\begin{equation}\label{ONSH0}
 H_\sigma f = \sum_{j=1}^{\infty}c_j\langle f,
 \varphi_j\rangle\varphi_j .
\end{equation}
It  follows that
\begin{equation}\label{ONSHlam}
H_\lambda f = H_{T_\lambda \sigma} f =  \pi(\lambda)H_\sigma \, \pi(\lambda )^{\ast} f = \sum_{j=1}^{\infty}c_j\langle
f,
  \pi(\lambda )\varphi_j\rangle \pi(\lambda )\varphi_j,
\end{equation}
and  $\{\pi (\lambda) \varphi_j, j\in \bN \}$ is an orthonormal system of eigenfunctions
of $H_\lambda$.
\end{Cor}

A priori, the spectral representation of $H_\lambda$ holds only for
$f\in\mathbf{L}^2 (\mathbb{R}^d )$. The next corollary extends  the
spectral representation to all of  $\mathbf{M}^\infty _{1/\nu} (\mathbb{R}^d
)$.
 
\begin{Cor}\label{Cor:SpecrepMinf}
The expansion for $H_\lambda f$ given in \eqref{ONSHlam} is
well-defined on  $ \mathbf{M}^\infty _{1/\nu } (\mathbb{R}^d)$ and converges to
$H_\lambda f$  in $\mathbf{L}^2$ for all $f\in
\mathbf{M}^\infty _{1/\nu }(\mathbb{R}^d)$. 
\end{Cor}
\begin{proof}
Without loss of generality, we assume $\lambda = 0$ and set $H =
H_\sigma $.
Since $Hf \in  \mathbf{L}^2(\rd ) $ for
every $f\in \mathbf{M}^\infty _{1/\nu}(\rd )$ by Lemma~\ref{Le:PropH_0}(iii), we can
expand $Hf$ with respect to the orthonormal system of eigenfunctions
of $H$ and obtain that 
\begin{equation}\label{Eq:serMin}
	Hf = \sum_{j=1}^\infty \langle Hf, \varphi_j\rangle \varphi_j
        + r 
\end{equation}
for some $r\in \mathbf{L}^2(\rd )$ in the orthogonal complement of
$\mathrm{span}\, \{ \vf _j: j\in \bN \}$. As $H$ is self-adjoint on
$\mathbf{L}^2(\rd )$, we also have $\langle Hf, \vf _j\rangle =
\langle f, H\vf _j\rangle = c_j \langle f, \vf _j\rangle $, and
consequently
\begin{equation}
  \label{eq:cc15}
  Hf =\sum_{j=1}^\infty c_j\langle f, \varphi_j\rangle     \varphi_j
  + r  \, .
\end{equation}
We need to show that  $r= 0$. Since $r\in \mathbf{L}^2(\rd )$ is
orthogonal to all eigenfunctions $\vf _j$,  we find that $\langle Hf,
r\rangle = \|r\|_2^2$. 

To show $r=0$, we first observe that  $\langle Hh,r\rangle = 0$ for
all $h\in \mathbf{L}^2(\rd )$ by \eqref{ONSH0}.  
Since $\mathbf{L}^2(\rd )$ is $w^*$-dense in $\mathbf{M}^\infty
_{1/\nu } (\rd
)$, we may choose  an approximating sequence $f_n\in \mathbf{L}^2(\rd
)  $  such that $f_n \stackrel{w^*}{\to} f\in \mathbf{M}^\infty
_{1/\nu }(\rd )$. For instance, $f_n$ may be chosen as
\[f_n = \int_\rdd \chi_{B_n}(z) \cV_g f(z) \pi (z) g\, dz, \] where $B_n$ is
the ball  with radius $n$ and  centered at $0$. 
 Furthermore, since
$f_n \stackrel{w^*}{\to} f$, we obtain in particular that 
$\cV_\varphi f_n$ converges to $  \cV_\varphi f$ uniformly on compact
sets~\cite[Theorem~4.1]{fg89jfa}. Consequently   
\[0= \langle Hf_n,r\rangle  =  \int_{\mathbb{R}^{2d}}\sigma (z)
\cV_\varphi f_n (z)\overline{\cV_\varphi r (z)}\, dz \,\, \to \int_{\mathbb{R}^{2d}}\sigma (z)
\cV_\varphi f (z)\overline{\cV_\varphi r (z)}\, dz = \langle Hf,
r\rangle = \|r\|_2^2 \, .\] 
This shows that $r=0$ and so 
the series  \eqref{ONSHlam} represents $Hf$ for all $ f\in
\mathbf{M}^\infty _{1/\nu } (\rd )$. 
 \end{proof}


\section{From Local Information to Global Information}\label{Se:Main}
We first  state and prove the  main result for the modulation spaces
$\mathbf{M}^p_m (\mathbb{R}^d )$. The  generalizations to $\MS $ 
will be discussed later. As always,  $\nu $ denotes a submultiplicative,
even weight function on $\rdd $ satisfying the condition $\lim _{n\to
  \infty } \nu (nz)^{1/n} =1$ for all $z\in \rdd $. 
  
\begin{Th}\label{Th:Main}
Let $\sigma\in \mathbf{L}^1 _\nu (\bR^{2d})$ be a non-negative symbol satisfying the condition 
\begin{equation}\label{eq:symbcond}
A\leq\sum_{\liL} T_\lambda \sigma\leq B,\,\,   a.e.
\end{equation}
for two constants $A,B>0$.  Assume that $\varphi\in \mathbf{M}^1_\nu
(\bR^d )$. 
Then for every  $\nu$-moderate weight $m$ and $1\leq p<\infty $ the
distribution  $f\in \mathbf{M}^\infty _{1/\nu }(\rd )$ belongs to $\mathbf{M}^{p}_m (\rd )$,   if and only if 
\begin{equation}\label{eq:Mpcond}
\Big( \sum_{\liL} \| H_\lambda f\|_2^p \, m(\lambda )^p \Big) ^{1/p} < \infty \, ,
\end{equation}
and the expression in  \eqref{eq:Mpcond} is  an equivalent norm on
$\mathbf{M}^p_m(\bR^d )$. 

Similarly, for  $p=\infty $   we obtain  the norm equivalence
\begin{equation}\label{eq:Minfcond}
\| f \|_{\mathbf{M}^{\infty}_m}\asymp\sup_{\liL} \| H_\lambda f\|_2 \,
m(\lambda ) .
\end{equation}
\end{Th}

The norm equivalence supports the interpretation that $H_\lambda f$
carries the local \tf\ information about $f$ near $\lambda \in \rdd
$. By combining the local pieces $H_\lambda f$, one obtains the global
\tf\ information as it is measured by modulation space norms. 


The proof of Theorem~\ref{Th:Main} requires some preparations.
We first show that  finitely many eigenfunctions of $H_0=\cV _\vf ^*
\sigma \cV _\vf $   generate a
multi-window Gabor frame for $\mathbf{L}^2(\rd )$. With this crucial
step in place, Theorem~\ref{Th:Main} can then  be deduced from the characterization
of \modsp s by means of Gabor frames. 

\subsection{Multi-Window Gabor Frames}
\begin{Lem}\label{Lem:Fineig} Assume that $\sigma \in
  \mathbf{L}^1(\rdd) $ and $\sum _{\lambda \in \Lambda } T_\lambda
  \sigma \asymp 1$,
  and that $\vf \in \mathbf{M}^1_\nu  (\rd )$.  
Let $\{\varphi_j : j\in\bN\}$ be the orthonormal system of
eigenfunctions of $H_0$. Then there exists $n\in \bN$, such that 
the finite union $\bigcup_{j=1}^n\mathcal{G} (\varphi_j,
\Lambda )$ is a multi-window Gabor frame for $\Lz$.  
\end{Lem}
An analogous  statement was proved and used in~\cite{dofegr06} for the
lattice $\Lambda = \zdd $ and rational lattices by means of Zak
transform methods. 
In the case of general lattices  we cannot apply Zak-transform
methods. As a substitute, we will
use  a  finite intersection property for $\Lambda $-invariant
subspaces of $\mathbf{M}^\infty$. The following statement may be  of interest
in its own right.

\begin{Lem}\label{Le:cpL}
Assume that $\mathcal{W}_n$ is a  sequence of $w^{\star}$-closed subspaces in $\mathbf{M}^{\infty} (\rd )$ such that
\begin{itemize}
\item[(i)] $\mathcal{W}_n \supseteq \mathcal{W}_{n+1}\neq\{0\}$ for all $n\in\mathbb{N}$ and 
\item[(ii)] $\mathcal{W}_n$ is invariant under all operators $\pi
  (\lambda)$ for  $\lambda \in \Lambda $.
\end{itemize}
Then $\bigcap_{n\geq 1} \mathcal{W}_n \neq \{0\}$.

\end{Lem}

\begin{proof}
Let $Q$ be the closure of  a relatively compact  fundamental domain of
$\Lambda$, for instance, if $\Lambda = A\zdd $, then
$Q=A[0,1]^{2d}$. We first  choose  a sequence
$h_n\in\cW_n$ with $\|h_n\|_{\mathbf{M}^{\infty}} = \sup _{z\in \rdd }
|\cV _\vf h_n(z)| =1$. 
Then there exists a sequence of points $\lambda_n$ in $\Lambda$,
such that  
\[\sup_{z\in Q} |\cV_{\varphi} (\pi(\lambda_n )h_n)(z)| = 1\, .\] 
Since $\mathcal{W}_n$ is invariant under all $\pi (\lambda), \lambda
\in \Lambda ,$ the distribution  $f_n = \pi (\lambda_n) h_n$ is in
$\mathcal{W}_n$. \\
Next  we show that the set of restrictions $\{\cV_{\varphi} f_n |_Q\}$
is equicontinuous. 
We have 
\begin{equation}
	|\cV_\varphi f_n (z)-\cV_\varphi f_n (\xi)| = |\langle f_n,
        (\pi (z)-\pi ( \xi ))\varphi \rangle|\leq
        \|f_n\|_{\mathbf{M}^\infty}\cdot\|(\pi (z)-\pi ( \xi
        ))\varphi\|_{\mathbf{M}^1} \, .
\end{equation}
Since $\|f_n \|_{\mathbf{M}^\infty } = \|\pi (\lambda
_n)h_n\|_{\mathbf{M}^\infty } = 1$,  the equicontinuity follows from the  strong continuity of
time-frequency shifts on $\mathbf{M}^1 (\mathbb{R}^d )$. \\ 
We next choose $z_n \in Q$ with $ |\cV_{\varphi} f_n(z_n )|\geq
\frac{1}{2}$. Since the unit ball in $\Min (\rd )$ is
$w^{\star}$-compact, there exists a subsequence $f_{n_k}$ that  converges
to some  $f\in \mathbf{M}^\infty (\rd )$ in the $w^{\star}$-sense.  Furthermore, by compactness of $Q$,
there also exists a subsequence $z_{\ell}$ of $z_{n_k}$, such that
$z_{\ell }\rightarrow z\in Q$. Hence, by equicontinuity,  
\[\cV_{\varphi} f_{\ell }(z_{\ell } )\rightarrow \cV_{\varphi} f(z )\, .\]
Since  $|\cV_{\varphi} f_{\ell }(z_{\ell } )|\geq 1/2$, we  conclude that
also $|\cV_{\varphi} f(z )|\geq 1/2$, and consequently $f\neq 0$. 

By construction, $f_{\ell }\in\cW_m$ for every  $\ell \geq m $, hence
 we obtain  $f= w^*-\lim _{\ell
  \to \infty } f_\ell \in\cW_m$  for all $m$, because  $\cW_m$ is $w^{\star}$-closed.
To summarize,  we have constructed a non-zero  $f\in\Minf$  that is in
$\cW_m$ for all $m$.  

\end{proof}

\begin{proof}[Proof of Lemma~\ref{Lem:Fineig}]
To prove that finitely many eigenfunctions generate a \mwgf\ with
respect to the lattice $\Lambda $, we assume on the contrary that
$\bigcup_{j=1}^n\mathcal{G} (\varphi_j, \Lambda )$ is not a frame 
  for every  $n\in\mathbb{N}$. Using Lemma~\ref{Le:cpL} and
  Lemma~\ref{Le:MuwiwithoutIneq}, we 
  will derive  a contradiction to the assumption that $A\leq \sum _{\lambda
    \in \Lambda } T_\lambda \sigma \leq B$. 

We use the criterion of Lemma~\ref{Le:MuwiwithoutIneq}. Let
$\boldsymbol\varphi_n  = (\varphi_1,\ldots,\varphi_n)$ be the 
  vector-valued function consisting of  the first $n$ eigenfunctions of 
  $H_0$, and 
\[\cW_n = \ker (C_{\boldsymbol\varphi_n,\Lambda}) = \{f \in\Minf:
\langle f,\pi (\lambda )\varphi_j\rangle = 0,  \ \forall
\lambda\in\Lambda , j = 1,\ldots , n\}\]  
be the kernel of the coefficient operator $C_{\boldsymbol\varphi_n
  ,\Lambda}$ in $\mathbf{M}^\infty (\rd )$. 

If $\bigcup_{j=1}^n\mathcal{G} (\varphi_j, \Lambda )$ is not a frame,
then  $\cW _n$ is a non-trivial subspace of $\mathbf{M}^\infty (\rd )$
by Lemma~\ref{Le:MuwiwithoutIneq}. By construction, the $\cW 
_n$'s form a nested sequence of w$^*$-closed subspaces of
$\mathbf{M}^\infty (\rd )$, and they are also invariant under $\pi
(\lambda ), \lambda \in \Lambda $. Thus the assumptions of
Lemma~\ref{Le:cpL} are satisfied, and we conclude that $\bigcap
_{n=1}^\infty \cW _n \neq \{0\}$.  This means that there exists a
non-zero $f\in \mathbf{M}^\infty (\rd )$, such that 
\begin{equation}
\label{Eq:zeros}
 \langle f,\pi (\lambda
  )\varphi_j\rangle = 0 \qquad  \mbox{  for all } \lambda \in \Lambda
  \mbox{ and all
  }j\in\mathbb{N} . 
\end{equation} 
We now consider $H_\lambda f$.  Since $H_\lambda f\in \mathbf{M}^1
(\mathbb{R}^d)$ by Lemma~\ref{Le:PropH_0}, the bracket $\langle
H_\lambda f, f\rangle$ is well-defined and  given by 
\begin{equation}\label{Eq_Hla1}
	\langle H_\lambda f, f\rangle = \int_{\mathbb{R}^{2d}} \sigma (z-\lambda ) |\cV_\varphi f (z ) |^2 dz .
\end{equation}
On the other hand, the extended spectral representation of  Lemma~\ref{Cor:SpecrepMinf} and \eqref{Eq:zeros}
imply that 
\begin{equation}\label{Eq:rep1}
H_\lambda f = \sum_{j = 1}^\infty c_j \langle f, \pi (\lambda )
\varphi _j \rangle \pi (\lambda )  \varphi _j  = 0. 
\end{equation}
Consequently $\langle H_\lambda f,f \rangle = 0$ for all $\lambda \in
\Lambda $, and   $|\cV_\varphi f (z ) |^2$ vanishes on  $\bigcup
_{\lambda \in \Lambda } \mathrm{supp}\, T_\lambda \sigma$.

According to the crucial  assumption \eqref{eq:symbcond} we have
$\sum_{\lambda\in\Lambda} T_\lambda \sigma \geq A >0$ almost
everywhere, and thus  $\bigcup_{\lambda\in\Lambda} \supp (T_\lambda \sigma) = \mathbb{R}^{2d}$. 
Therefore, \eqref{Eq_Hla1} and \eqref{Eq:rep1} imply that  $\cV_\varphi f =
0$, from which $f = 0$ follows. 
This is a contradiction to $f$ being a non-zero element in $\bigcap
_{n=1}^\infty \cW _n$.

This contradiction shows that there exists an $n\in \bN $, such that $
\bigcup_{j=1}^n\mathcal{G} (\varphi_j, \Lambda )   $  is  a multi-window Gabor frame, and we are done. 
 \end{proof}
\begin{rem}
Note that for finite-rank operators $H_0$, it can be seen directly
that the finite set of eigenvectors generates a multi-window Gabor
frame for $\Lambda$.  
\end{rem}
\subsection{Proof of Theorem~\ref{Th:Main} }
We are now ready to prove the main theorem.
We observe that for $f\in \bm _{1/\nu }^\infty (\rd )$, 
$H_\lambda f \in \mathbf{L}^2 (\rd ) $  by
Lemma~\ref{Le:PropH_0}(iii). Thus the terms in \eqref{eq:Mpcond} are 
well-defined. 

 First assume  that $p<\infty$ and  $f\in \mathbf{M}^p_m(\mathbb{R}^d )\subseteq
 \bm ^\infty _{1/\nu 
}(\rd )$. Using the embedding $ \mathbf{M}^1 (\mathbb{R}^d
 )\hookrightarrow\Lz$ and  the estimate \eqref{regineq} with $\nu
 \equiv 1$, we majorize   $\| H_{\lambda} f\|_2 $ as follows:
\begin{align}\label{Lpest1}
\| H_\lambda f\|_2  &\leq C_{\varphi}\| H_\lambda f\|_{\mathbf{M}^1} \\\notag
   & \leq C_{\varphi}
\|(T_\lambda\sigma ) \cdot \mathcal{V}_{\varphi}f\|_1 \, \|\mathcal{V}_{\varphi}\varphi\|_1 \\\notag
& =   C_{\varphi} C\int_\rdd   |\sigma(z-\lambda)|\cdot
|\mathcal{V}_{\varphi}f(z)|\, dz\\\notag 
& =  C_{\varphi}C (|\mathcal{V}_{\varphi}f|\ast \sigma^\vee )(\lambda ),\notag
\end{align}
where $\sigma^\vee (z ) = \sigma^\vee (-z)$.
Thus  $\| H_{\lambda} f\|_2 $ is majorized by a sample of
$|\mathcal{V}_{\varphi}f|\ast \sigma^\vee$. To proceed further, we
use  the fact that $\mathcal{V}_{\varphi}f\in  \Wsp (\mathbf{L}^p_m)$
and $\|\mathcal{V}_{\varphi}f\|_{\Wsp (\mathbf{L}^p_m)}\leq C_0
\|\vf \|_{\mathbf{M}_\nu ^1} \|f\|_{\mathbf{M}^p_m}$ for  $\varphi\in
\mathbf{M}^1_\nu (\rd )$ and $f\in \mathbf{M}^p_m 
(\rd )$ by~\cite[Thm.~12.2.1]{gr01}.  Now  the convolution relation
\eqref{Eq:WAsamp} and the sampling inequality \eqref{Eq:lpsamp0} imply
that   
\begin{align}\label{Lpest2}
 \sum_\lambda\| H_\lambda f\|^p_2 \,  m(\lambda )^p &\leq  C_{\varphi}C  \|(\sigma^\vee \ast
|\mathcal{V}_{\varphi}f|)|_\Lambda\|^p_{\ell^p_{{m}}}\\\notag
&\leq  C_{\varphi}C  C_\Lambda \|\sigma^\vee \ast
|\mathcal{V}_{\varphi}f|\|^p_{\Wsp (\mathbf{L}^p_m)}\\\notag
&\leq  C_{\varphi}C  C_\Lambda 
\|\sigma\|^p_{\mathbf{\mathbf{L}^1_\nu}} \|\mathcal{V}_{\varphi}f\|^p_{\Wsp (\mathbf{L}^p_m)}\leq C_{\varphi}C C_\Lambda 
\|\sigma\|^p_{\mathbf{\mathbf{L}^1_\nu}} \|f\|^p_{\mathbf{M}^p_m}.
\end{align}
The same argument yields $\sup _{\lambda \in \Lambda } \|H_\lambda
f\|_2 \, m(\lambda ) \leq C \|f\|_{\mathbf{M}^\infty _m}$. 

 Hence, for $1\leq p \leq \infty$, the mapping $f\rightarrow
 (\|H_\lambda f\|_2)_{\lambda\in\Lambda}$ is bounded from
 $\mathbf{M}^p_m (\rd )$ to $\ell^p_{{m}} (\Lambda )$. 

Conversely, assume  that  $p<\infty$ and
\[ \sum_\lambda\| H_\lambda f\|^p_2 m(\lambda )^p <\infty.\] 
We  need to show that $f\in \mathbf{M}^p_m (\rd )$.
Since $\|H_\lambda f\|_2  = \sup_{\|g\|_2 = 1}|\langle H_\lambda f,g\rangle |$, we have the inequality
\[\sum_{\lambda}|\langle H_\lambda f,g_\lambda \rangle |^p
m(\lambda)^p\leq \sum_{\lambda}\|H_\lambda f\|^p_2 \,  m (\lambda)^p
<\infty\] 
for arbitrary sequences $g_\lambda\in\mathbf{L}^2 (\rd )$ with
$\|g_\lambda\|_2  = 1$.  Applying the eigenfunction expansion of
Corollary~\ref{Cor:Specrep}, we obtain 
\begin{equation}\label{Eq:estj}
\sum_{\lambda}\left|\sum_{j=1}^\infty c_j\langle f,\pi (\lambda
  )\varphi_j\rangle \langle \pi (\lambda )\varphi_j
  ,g_\lambda\rangle\right|^p m (\lambda)^p\leq
\sum_{\lambda}\|H_\lambda f\|^p_2 m (\lambda)^p<\infty .\end{equation} 
Now fix $j_0\in\mathbb{N}$ and set $g_\lambda = \pi (\lambda
)\varphi_{j_0}$ for $\liL$. Since the eigenfunctions of $H_\lambda$
are orthonormal, the sum over $j$ collapses to a single term,  and \eqref{Eq:estj}
becomes 
\begin{equation*}\label{Eq:lj0est}
\sum_{\lambda}|\langle H_\lambda f,g_\lambda \rangle |^p m (\lambda)^p = \sum_{\lambda}|c_{j_0}\langle  f,\pi (\lambda )\varphi_{j_0}\rangle |^p m (\lambda)^p\leq  \sum_\lambda\| H_\lambda f\|^p_2 m(\lambda )^p <\infty . 
\end{equation*}
 The last inequality holds for every $j_0\in \bN$.
After summing over finitely many $j_0$ and switching to the  $\ell ^2$-norm
on $\bC ^n$, we obtain the inequality 
\begin{eqnarray} \notag 
\lefteqn{\sum_{\lambda}   \Big(   \sum_{j=1}^n|\langle  f,\pi (\lambda
)\varphi_j\rangle |^2\Big)^{1/2} m (\lambda) ^p \leq 
  \sum_{j=1}^n\sum_{\lambda}|\langle  f,\pi (\lambda
)\varphi_j\rangle |^p m (\lambda)^p }   \\
&\leq & \Big(\sum_{j=1}^n
\frac{1}{c_j^p}\Big) \, \sum_\lambda\| H_\lambda f\|^p_2 m(\lambda )^p
<\infty.  \label{eq:c18}
\end{eqnarray}
We now apply  Lemma~\ref{Lem:Fineig} and choose  an $n\in \bN$,   such that 
 $\bigcup_{j = 1}^n\mathcal{G} (\varphi_j ,\Lambda ) $ is a
 multi-window Gabor frame for $\mathbf{L}^2(\mathbb{R}^d )$. Since all
 $\vf _j$ are in $\mathbf{M}_\nu ^1(\rd )$, the fundamental
 characterization of \modsp s (Section~\ref{Se:charMSMW}) is valid. Thus
 Theorem~\ref{fund}(i) implies that $f\in \mathbf{M}^p_m(\rd )$.

If  $p = \infty$ and   $\sup_{\liL} \| H_\lambda f\|_2\,  m(\lambda
)<\infty$, then, by choosing $g_\lambda$ as before, we find 
\[c_{j_0} \sup_\lambda |\langle f,\pi (\lambda ) \varphi_{j_0}\rangle
|m(\lambda )\leq \sup_\lambda \|H_\lambda f\|_2 m(\lambda )<\infty\] 
for every $j_0$.

Arguing as above,  Theorem~\ref{fund} says that 
\[\|f\|_{\mathbf{M}^\infty_m}\leq C\, \max_{j=1,\ldots ,
  n}\sup_\lambda |\langle f,\pi (\lambda ) \varphi_{j}\rangle
|m(\lambda )\leq \Big(\max_{j=1,\ldots , n}\frac{1}{c_j}\Big) \sup_\lambda
\|H_\lambda f\|_2 m(\lambda )<\infty,\]
and $f\in \mathbf{M}^\infty _m(\rdd )$. 

Combining \eqref{Lpest2} and \eqref{eq:c18}, we have shown that  
 $\|f\|_{\mathbf{M}^p_m}$ and $\Big(\sum_{\liL} \| H_\lambda f\|_2^p
 \, m(\lambda )^p\Big)^{1/p} $  for $1\leq p<\infty $ (or $\sup_{\liL} \| H_\lambda f\|_2
 m(\lambda )$ for $p=\infty$)  are equivalent norms on $\mathbf{M}^p_m
 (\rd )$.  \hfill \qed 
 \subsection{Variations of  Theorem~\ref{Th:Main}} 
In order to formulate our main result for mixed-norm spaces and
arbitrary lattices, we have to resort to  the theory of coorbit
spaces, as introduced in~\cite{fg89jfa,fg89mh}. In particular, for
arbitrary lattices, a sequence $(c_\lambda)_{\lambda\in\Lambda}$ is in
the sequence spaces associated with $\mathbf{L}^{p,q}_m (\rdd )$, if
$\sum_{\lambda\in\Lambda} c_\lambda\chi_{\lambda+ Q}$ is in
$\mathbf{L}^{p,q}_m (\rdd )$ for some  fundamental domain $Q$  of
$\Lambda$.  
With this definition, we may give the following characterization.  
\begin{Th}\label{Cor:coovers}
Let $\Lambda$ be an arbitrary lattice in $\rdd$  and $Q$ be  a
relatively compact   fundamental domain $Q$. Assume  the same
conditions  on  $\sigma$ and $\varphi$ as in Theorem~\ref{Th:Main}. 
Then a distribution  $f\in \bm ^\infty _{1/\nu }(\rd )$ belongs to $
\mathbf{M}^{p,q}_m (\rd )$, $1\leq p,q \leq\infty$,  if and only if  
\begin{equation}\label{eq:Mpcondcoorb}
\sum_{\liL} \| H_\lambda f\|_2\chi_{\lambda+Q}\in\mathbf{L}^{p,q}_m(
\mathbb{R}^{2d}) \, ,
\end{equation}
and  $\big\|\sum_{\liL} \| H_\lambda f\|_2\, \chi_{\lambda+Q}\big\|_{\mathbf{L}^{p,q}_m} \asymp \|f\|_{\mathbf{M}^{p,q}_m}$.  
\end{Th}

\begin{proof}
  The proof is almost identical to the proof of
  Theorem~\ref{Th:Main}. The only modifications occur in \eqref{Lpest2},
  which has to be replaced by 
$$
\|\sum _{\lambda \in \Lambda } \|H_\lambda f\|_2 \, \chi _{\lambda +Q}
\|_{\Lmpq } \leq \|\sum _{\lambda \in \Lambda } |\cV _\vf f \ast
\check{\sigma }(\lambda )| \, \chi _{\lambda +Q} \|_{\Lmpq } \leq C \|\cV _\vf f \ast
\check{\sigma } \|_{\mathbf{W}(\Lmpq )} \, .
$$
Likewise, in \eqref{eq:c18}  we replace the weighted
$\mathbf{L}^p_m$-norm by the general $\Lmpq $-norm. 
\end{proof}
For a  separable lattice $\Lambda = a\mathbb{Z}^d\times
b\mathbb{Z}^d$ the norm in \eqref{eq:Mpcondcoorb} is just the
$\ell ^{p,q}_{\tilde{m}}$-norm on $\zdd $ with $\tilde m (k,n) =
m(ak,bn)$.  In this case, $\lambda = 
(ka, nb)$, $k,n\in\mathbb{Z}^d$ and we may  write $H_\lambda f = H_{k,n}f$.
\begin{Cor}\label{Cor:seplatt}
Let $\Lambda = a\mathbb{Z}^d\times
b \mathbb{Z}^d$ be a separable lattice and assume  the  same
conditions on $\sigma$ and $\varphi$  as in Theorem~\ref{Th:Main}. 
Then a distribution $f\in \mathbf{M}^\infty _{1/\nu } (\rd )$ belongs to $ \mathbf{M}^{p,q}_m (\rd )$ for $1\leq p,q<\infty$,  if and only if 
\begin{equation}\label{eq:Mpcondseplatt}
 \Big(\sum_{n\in\mathbb{Z}^d}   (\sum_{k\in\mathbb{Z}^d}\|H_{k,n}f\|^p_2
 \, m(ka,nb)^p)^{q/p}\Big)^{1/q} <\infty \, ,
\end{equation} 
and  \eqref{eq:Mpcondseplatt} defines an equivalent norm on
$\mathbf{M}^{p,q}_m(\mathbb{R}^d )$. The result holds for $p=\infty$
or $q =\infty$ with the usual  modifications.  
\end{Cor}

\subsection{Existence of multi-window Gabor frames and properties of the eigenfunctions $\varphi_j$}\label{Se:Con}
We finally point out some immediate consequences of our results and
methods.

%

The intermediate results leading to Theorem~\ref{Th:Main} also imply
the existence of \mwgf s for general lattices. 

\begin{Th}\label{Pro:FinGenLatt}
Let $\Lambda$ be an arbitrary lattice and $\nu$ a submultiplicative
weight on $\mathbb{R}^{2d}$. Then there exist finitely many functions
$\varphi_j\in \mathbf{M}^1_\nu (\mathbb{R}^d)$, such that
$\bigcup_{j=1}^n\mathcal{G} (\varphi_j, \Lambda )$ is a multi-window
Gabor frame for $\mathbf{L}^2 (\mathbb{R}^d )$. 
\end{Th}
\begin{proof}
Choose $\sigma\in\mathbf{L}^1_\nu (\mathbb{R}^{2d} )$ such that
$\sum_{\liL} T_\lambda \sigma\asymp 1$ and fix a window $\vf \in
\mathbf{M}^1_\nu (\rd )$.  For instance, one may choose
the characteristic function $\chi _Q$ of a (relatively compact)
fundamental domain of $\Lambda $ and the Gaussian window $\vf (t) =
e^{-\pi t\cdot t}$. 

Now   consider  the localization operator 
$H_0 = \cV_\varphi^\ast \sigma \cV_\varphi$. 
According to Lemma~\ref{Le:PropH_0}(ii),  all  eigenfunctions
$\varphi_j$ of $H_0$  belong to $\mathbf{M}^1_\nu (\mathbb{R}^d )$.
Lemma~\ref{Lem:Fineig} states that for some finite $n\in \bN$ the set
$\bigcup_{j = 1}^n\mathcal{G} (\varphi_j ,\Lambda ) $ is a 
 multi-window Gabor frame for $\mathbf{L}^2(\mathbb{R}^d )$. 
\end{proof}

The existence of \mwgf s for general lattices was known before. 
On the one hand, it is an immediate consequence of coorbit theory
applied to the Heisenberg group. 
To be more precise, according to   \cite[Thm.~7]{fg92chui}
 for every  lattice $\Lambda $ and every non-zero   $g\in
 \mathbf{M}^1_\nu (\rd )$ there exists $n\in \bN $, such that  the set $\cG (g,
\frac{1}{n}\Lambda )$ is a Gabor frame for $\lrd $. Using a coset
decomposition  $\frac{1}{n}\Lambda = \bigcup (\mu +\Lambda )$ for
suitable $\mu \in \Lambda $, one 
sees that     $\cG
(g, \frac{1}{n}\Lambda ) = \bigcup \cG (\pi (\mu )g, \Lambda )$ is a
\mwgf\ with all windows $\pi (\mu )g$ derived from a single window $g$. 
 Recently Luef~\cite{lu09} proved  the existence of  
 multi-window Gabor
 frames by  exploiting a  connection between Gabor analysis and  non-commutative
 geometry. Our methods  provide  a third,  independent proof
 for this interesting result.\\

 The construction of
 \mwgf s in Proposition~\ref{Pro:FinGenLatt} yields more detailed
information about the  frame  generators, since   they are eigenfunctions of a localization operator.
 Intuitively the  eigenfunctions corresponding to the
largest eigenvalues of a  localization operator  concentrate their
energy on the essential support of the symbol $\sigma $ of $H_0$.  For
the special case of compactly supported   $\sigma  $, 
this intuition is made  precise by the following result. 
 
\begin{Pro}\label{Pro:EigCon}Let the  non-negative function
  $\sigma\in\mathbf{L}^1(\mathbb{R}^{2d})$ be supported  in  
a
compact set $\Omega$ in $\mathbb{R}^{2d}$ with $0\leq\sigma(z) \leq C_\sigma<\infty$ for $z\in\Omega$.
 Consider the localization operator  given by $H_\sigma f = \cV_\varphi 
^\ast \sigma \cV_\varphi f$ with  $\vf \in \mathbf{M}^1(\rd
)$,   $\|\varphi\|_2=1$ and spectral representation  as in
Corollary~\ref{Cor:Specrep}. Then  the eigenfunctions $\vf _j$  of
$H_\sigma $    satisfy the following time-frequency concentration
\begin{equation}
 \int_\Omega |\cV_\varphi \varphi_j  (z)|^2\, dz\geq \frac{c_j}{C_\sigma}\, .
\end{equation}

Equality holds, if and only if  $\sigma(z) /C_\sigma = \chi_\Omega(z)$ is the characteristic  function of $\Omega$.\end{Pro}
\begin{proof} Using the weak interpretation of $H_\sigma $ from~\eqref{eq:f1}, we obtain 
\begin{align*}
	\int_{\Omega}|\cV_\varphi \varphi_j  (z)|^2\,
        dz\geq&\frac{1}{C_\sigma}\int_{\Omega}\sigma (z)|\cV_\varphi
        \varphi_j(z)|^2 \, dz\\	  
=&\frac{1}{C_\sigma}\langle H_\sigma \varphi_j, \varphi_j\rangle
 =\frac{c_j}{C_\sigma}\| \varphi_j\|_2^2=\frac{c_j}{C_\sigma}\, .
\end{align*}
\end{proof}

\begin{appendix}
  \section{Characterizations of  Modulation Spaces and Multi-Window Gabor Frames}\label{App}
In the  appendix, we will  sketch  the  proof of
Theorem~\ref{fund} and formulate a series of new characterizations of
multi-window Gabor frames. These statements 
generalize  well-known facts from Gabor analysis  and the
results about Gabor frames without inequalities in~\cite{gr07}. 

For the investigation of multi-window Gabor frames we need  the dual
concept of vector-valued Gabor systems.  
 In this case we consider the Hilbert space $\mathcal{H} = \Lt (\bR^d,\bC^n )$ consisting of all vector-valued functions $\mathbf{f} (t) = (f_1(t),\ldots ,f_n (t) )$ with the inner product
 
\begin{equation}
	\langle \mathbf{f} , \boldsymbol\varphi\rangle
        _{\mathbf{L}^2(\rd , \bC ^n)} = \sum_{j = 1}^n \int f_j
        (t)\overline{ \varphi_j (t)} dt = \sum_{j = 1}^n\langle f_j,
        \varphi_j\rangle _{\mathbf{L}^2(\rd)}. 
\end{equation}
Time-frequency-shifts act coordinate-wise on $\mathbf{f}$. The
vector-valued Gabor system $\mathcal{G}(\boldsymbol\varphi,\Lambda ) =
\{ \pi (\lambda )\boldsymbol{\varphi} : \lambda \in \Lambda \}$ 
is a Riesz sequence in  $\Lt (\bR^d, \bC^n  )$, if there exist
constants $0<A,B<\infty$ such that   
for  all finitely supported sequences $\mathbf{c}$,
\begin{equation}\label{Eq:RieszCond}
	A\|\mathbf{c}\|_2^2\leq \|\sum_{\lambda \in\Lambda} c_\mu \pi
        (\lambda  ) \g\|^2_{\mathbf{L}^2 (\mathbb{R}^d,\mathbb{C}^n)}\leq B\|\mathbf{c}\|_2^2.
\end{equation}

We now proceed to the  proof of Theorem~\ref{fund}. The crucial step
is to show  the invertibility of the frame operator on $\mathbf{M}^1_\nu 
(\mathbb{R}^{d})$. This step requires a special representation of the
frame operator due to Janssen~\cite{janssen95}  and at its core uses
``Wiener's lemma for twisted convolution''~\cite{grle04}. 

For  $\varphi_j, \phi_j$ in $\mathbf{M}^1(\bR^d), j=1, \dots , n,$ we
denote   frame-type operators by 
\[S_{\g,{\boldsymbol\psi}} f  =  \sum_{\liL}\sum_{j=1}^n\langle f,
\pi (\lambda )\varphi_j\rangle  \pi (\lambda )\psi_j = \sum _{j=1}^n
S_{\vf _j, \psi _j}\, .
\]
The frame operator of the Gabor system $\bigcup _{j=1}^n \cG (\vf _j,
\Lambda )$ is $S =S_{\g , \g}$. We usually omit the reference to the lattice
$\Lambda $ and the windows $\vf _j$. 

The  volume $s(\Lambda )$ of a lattice $\Lambda = A\zdd  $ is  defined as the measure
of a  fundamental domain of $\Lambda$ and is  $|\det
(A)|$.  
The adjoint lattice of $\Lambda $ is $\Lambda^\circ =
\{\mu\in\bR^{2d} : \pil\pi (\mu ) = \pi (\mu ) \pil\mbox{ for all }
\lambda\in \Lambda\}$.

\begin{Lem}[Janssen's representation]\label{Le:JannRep}
Assume that  $\varphi_j, \psi_j\in \mathbf{M}^1(\bR^d)$ for all
$j=1,\dots ,n$. Then    
the frame type operator associated to 
 $\bigcup_{j=1}^n\mathcal{G} (\varphi_j, \Lambda )$ and
 $\bigcup_{j=1}^n\mathcal{G} (\psi_j, \Lambda )$ 
can be written as 
\begin{equation}\label{Eq:MuWijann}
S_{\g,{\boldsymbol\psi}}f = s(\Lambda )^{-1}
  \sum_{\mu\in\Lambda^{\circ}}\sum_{j=1}^n\langle \varphi_j ,\pi (\mu
  )\psi_j\rangle \pi (\mu )f 
\end{equation}
with unconditional  convergence in the operator norm on $\mathbf{L}^2$.
\end{Lem}
\begin{proof}
By  Janssen's result~\cite{janssen95} the representation holds for a
single $S_{\vf _j, \psi _j}$ and \eqref{Eq:MuWijann} follows by taking
a sum. 
\end{proof}
The
\emph{canonical} dual frame is defined to be  $\gamma_{j,\lambda} =
\pil S^{-1}\varphi_j$  Since the frame operator $S= S_{\g , \g}$ commutes with
time-frequency shifts on  $\Lambda$,  we obtain  the reconstruction formulas
\begin{align*}
f =& S^{-1}S f = \sum_{\liL} \sum_{j =1}^n\langle f , \pi (\lambda
)\varphi_j\rangle \pil\gamma_j \\
=& S S^{-1} f
=\sum_{\liL} \sum_{j =1}^n\langle f , \pil \gamma_j\rangle \pil\varphi_j \\
=&D_{\boldsymbol\varphi,\Lambda}C_{\boldsymbol\gamma,\Lambda}f  = D_{\boldsymbol\gamma,\Lambda}C_{\boldsymbol\varphi,\Lambda}f
\end{align*}
As a general principle the localization of a frame is inherited by the
dual frame~\cite{FoG05}. The following statement is a generalization of
\cite[Thm.~9]{grle04} to multi-window Gabor frames on general
lattices.   

\begin{Lem}\label{Le:MuwiDual} Assume that $\nu $ is a
  submultiplicative, even weight on $\rdd $ satisfying $\lim _{n\to
    \infty } \nu (nz)^{1/n} =1 $ for all $z\in \rdd $. Assume further
  that   $\bigcup_{j=1}^n\mathcal{G} (\varphi_j, \Lambda )$ is a frame
  for $\mathbf{L}^2 (\mathbb{R}^d )$ and that $\varphi_j\in
  \mathbf{M}^{1}_\nu (\mathbb{R}^d )$. Then the frame operator $S$ is
  invertible on  $\mathbf{M}^{1}_\nu(\mathbb{R}^d )$ and $\gamma_j =
  S^{-1} \varphi_j\in \mathbf{M}^{1}_\nu (\mathbb{R}^d )$ for 
  $j=1, \dots , n$. 
\end{Lem}
\begin{proof}
Janssen's representation \eqref{Eq:MuWijann} implies that 
\begin{equation}
S = S_{\g, \g} = s(\Lambda )^{-1} \sum_{\mu\in\Lambda^{\circ}} c_\mu \pi (\mu ), 
\end{equation}
with a coefficient sequence $c_\mu  = \sum_{j=1}^n\langle \varphi_j
,\pi (\mu )\varphi_j\rangle$. The hypothesis     $\varphi_j\in
\mathbf{M}^{1}_\nu(\mathbb{R}^d )$ guarantees  that 
$\sum_{\mu\in\Lambda^{\circ}}|\langle \varphi_j, \pi (\mu )
\varphi_j\rangle | \nu (\mu ) <\infty$ for each $j$,
see~\cite[Cor.12.1.12]{gr01}, and therefore the coefficient sequence
$(c_\mu )$ is  in $\ell^1_{{\nu}}(\Lambda ^\circ )$.  Since 
$\bigcup_{j=1}^n\mathcal{G} (\varphi_j, \Lambda )$ is a frame, the
frame operator $S_{\g , \g}$ is invertible on $\mathbf{L}^2(\rd )$. It follows
from~\cite[Theorem~3.1]{grle04}  that the inverse frame operator
$S^{-1}$ is again  of the form $S^{-1} =
\sum_{\mu\in\Lambda^{\circ}} d_\mu \pi (\mu )$ with  a
coefficient sequence $\mathbf{d}$ in $\ell^1_{{\nu}}(\Lambda ^\circ )$.
This representation implies that $S\inv $ is bounded on
$\mathbf{M}^1_\nu (\rd )$ and that 
\begin{equation}
\|\gamma_j\|_{\mathbf{M}^1_\nu} =
\|S^{-1}\varphi_j\|_{\mathbf{M}^1_\nu} \leq C
\|\varphi_j\|_{\mathbf{M}^1_\nu} \, .
\end{equation} 
Therefore the  dual windows  $\gamma_j$, $j = 1,\ldots , n$ are in
$\mathbf{M}^1_\nu (\mathbb{R}^d)$ as claimed. 
\end{proof}
Once the  invertibility of the multi-window frame operator on
$\mathbf{M}^{1}_\nu(\mathbb{R}^d )$ is established, the proof of
Theorem~\ref{fund} is straight-forward by using the following
boundedness properties of the coefficient
 operator $C_{\g , \Lambda }$ and $D_{\g , \Lambda }$ from
 ~\cite[Theorem~12.2.3. and 12.3.4.]{gr01}. 
If   $\varphi_j \in \mathbf{M}^1_\nu(\mathbb{R }^d)$  and  $\gamma_j
\in \mathbf{M}^1_\nu(\mathbb{R }^d)$, then  both
$C_{\boldsymbol\varphi,\Lambda}$ and  $C_{\boldsymbol\gamma,\Lambda}$
are bounded from $\mathbf{M}^{p,q}_m (\mathbb{R }^d)$ into
$\ell_{{m}}^{p,q}(\Lambda, \mathbb{C}^n )$ for $1\leq p,q\leq \infty $
and for every  $\nu$-moderate
weight $m$. Likewise $D_{\boldsymbol\varphi,\Lambda}$ and 
$D_{\boldsymbol\gamma,\Lambda}$ are bounded from
$\ell_{{m}}^{p,q}(\Lambda, \mathbb{C}^n )$ into $\mathbf{M}^{p,q}_m
(\mathbb{R }^d)$. For the $\ell_{{m}}^{p,q}(\Lambda, \mathbb{C}^n
)$-norm we use the Euclidean norm on $\bC ^n$, so that
$\|\mathbf{c}\|_{\ell_{{m}}^{p,q}(\Lambda, \mathbb{C}^n )} = \|\sum
_{\lambda \in \Lambda } \big( \sum _{j=1}^n |c_{\lambda , j}|^2 \big)
^{1/2} \chi _{\lambda +Q} \|_{\mathbf{L}^{p,q}_m}$. 
 
As a consequence, the reconstruction formula  $f =
 D_{\boldsymbol\varphi,\Lambda}C_{\boldsymbol\gamma,\Lambda}f  =
 D_{\boldsymbol\gamma,\Lambda}C_{\boldsymbol\varphi,\Lambda}f$ holds
 for  $f\in\mathbf{M}^{p,q}_m (\mathbb{R}^d )$ with the correct norm estimates. The
 norm equivalence stated in Theorem~\ref{fund} then follows from   
\begin{align*}
\| f\|_{\MS}&= \|
D_{\boldsymbol\gamma,\Lambda}C_{\boldsymbol\varphi,\Lambda}f\|_{\MS}\leq
\|D_{\boldsymbol\gamma,\Lambda}\|_{op}\|
C_{\boldsymbol\varphi,\Lambda}f\|_{\ell^{p,q}_{{m}}(\Lambda , \bC ^n )}\\ 
&\leq\|D_{\boldsymbol\gamma,\Lambda}\|_{op}\|C_{\boldsymbol\varphi,\Lambda}\|_{op} \|f\|_{\MS}.
\end{align*}
\vspace{0.2cm}

 Next we come to  the characterization of \mwgf s
 (Lemma~\ref{Le:MuwiwithoutIneq})  and  extend the list of equivalent
 conditions.
For the formulation of the dual conditions on the adjoint lattice
$\Lambda ^\circ $ we  need the vector-valued versions of the analysis
and synthesis operators. For $\mathbf{f} = (f_1, \dots , f_n) \in
\mathbf{M}^\infty (\rd , \bC ^n)$ and $\g = (\vf _1, \dots , \vf _n)
\in \mathbf{M}^1(\rd, \bC ^n)$ the coefficient operator is defined to
be 
  $\widetilde C_{\boldsymbol\varphi,\Lambda^\circ}(\mathbf{f})(\mu)  =
     (\langle \mathbf{f},\pi (\mu ) \boldsymbol\varphi\rangle ), \mu
     \in \Lambda ^\circ$, and
     the synthesis operators is  
      $\widetilde D_{\boldsymbol\varphi,\Lambda^\circ}(\mathbf{c} ) =
      \sum_{\mu\in\Lambda^\circ}c_{\mu } \pi (\mu )
      \boldsymbol\varphi $.  
The Gramian operator 
$G_{\boldsymbol\varphi,\Lambda^\circ} =
\widetilde C_{\boldsymbol\varphi,\Lambda^\circ}\widetilde D_{\boldsymbol\varphi,\Lambda^\circ}$
 is  defined on  sequences indexed by
$\Lambda^\circ $.

  \begin{Lem}\label{Le:MuwiwithoutIneq1}Assume that  $\varphi_j\in \mathbf{M}^1 (\rd )$ for  $j = 1,\ldots , n$. The following are equivalent for the  multi-window Gabor system $\bigcup_{j=1}^n\mathcal{G} (\varphi_j, \Lambda )$:
\begin{itemize}
 \item[(i)] $\bigcup_{j=1}^n\mathcal{G} (\varphi_j, \Lambda )$ is a frame for $\Lz$.
\item[(ii)]Wexler-Raz biorthogonality: There exist  $\gamma_j\in \mathbf{M}^1(\mathbb{R}^d
  ), j=1,\dots, n$, such that
\begin{equation}\label{Eq:MuWiWexler}
s(\Lambda )^{-1} \sum_{j=1}^n\langle \varphi_j ,\pi (\mu )\gamma_j\rangle =\delta_{\mu ,0}\mbox{ for } \mu \in \Lambda^{\circ}.
\end{equation}
\item[(iii)]
Ron-Shen duality:  $\mathcal{G} (\boldsymbol\varphi, \Lambda^\circ )$
is a Riesz sequence  in $\mathbf{L}^2 (\mathbb{R}^d,\mathbb{C}^n)$.
  \item[(iv)] 
$S_{\g, \g}$    is invertible on $\mathbf{M}^1 (\bR^d )$.
   \item[(v)]$S_{\g, \g}$ 
     is invertible on $\mathbf{M}^{\infty} (\bR^d )$.
    \item[(vi)] $S_{\g, \g}$ 
      is one-to-one on $\mathbf{M}^{\infty} (\bR^d )$.
     \item[(vii)] The analysis operator $C_{\boldsymbol\varphi,\Lambda}: \mathbf{M}^{\infty} (\bR^d )\mapsto \ell^{\infty } (\Lambda , \bC^n )$ is one-to-one from $\mathbf{M}^{\infty} (\bR^d )$ to $\ell^{\infty } (\Lambda , \bC^n )$.
      \item[(viii)] The synthesis operator  $D_{\boldsymbol\varphi,\Lambda}$ defined on $\ell^1 (\Lambda ,\bC^n )$  has dense range in $\mathbf{M}^1 (\bR^d )$.
       \item[(ix)]$D_{\boldsymbol\varphi,\Lambda}$ is surjective from  $\ell^1 (\Lambda ,\bC^n )$ onto $\mathbf{M}^1 (\bR^d )$. 
        \item[(x)] The synthesis operator
          $\widetilde D_{\boldsymbol\varphi,\Lambda^\circ}$ defined on
          $\ell^\infty (\Lambda^\circ)$ 
      is one-to-one from  $\ell^\infty (\Lambda^\circ)$ to $\mathbf{M}^\infty (\bR^d,\bC^n )$.
         \item[(xi)] The analysis operator
           $\widetilde{C}_{\boldsymbol\varphi,\Lambda^\circ}$ defined on $\mathbf{M}^1(\bR^d, \bC^n )$ has dense range in  $\ell^1 (\Lambda^\circ )$ .        
\item[(xii)] $\widetilde C_{\boldsymbol\varphi,\Lambda^\circ}$ is surjective from  $\mathbf{M}^1(\bR^d, \bC^n )$ onto $\ell^1 (\Lambda^\circ)$.
\item[(xiii)] $G_{\boldsymbol\varphi,\Lambda^\circ}$ is invertible on $\ell^1(\Lambda^\circ )$.
\item[(xiv)]$G_{\boldsymbol\varphi,\Lambda^\circ}$ is invertible on $\ell^\infty(\Lambda^\circ )$.
\item[(xv)]$G_{\boldsymbol\varphi,\Lambda^\circ}$ is one-to-one on $\ell^1(\Lambda^\circ )$.
 \end{itemize}
\end{Lem}
The equivalence $(i) \Leftrightarrow (vi)$ is claimed in 
Lemma~\ref{Le:MuwiwithoutIneq} and is all we need for the main results
of our paper.  

\begin{proof}
The implication $(i)\Rightarrow (iv) $ was sketched  in Lemma~\ref{Le:MuwiDual}.  \\
$(i)\Leftrightarrow (ii) $:
 Time-frequency shifts on a lattice are linearly independent
in the following sense: if $c = (c_\mu
)_{\mu\in\Lambda^\circ}\in\ell^\infty$ and  
$\sum_{\mu\in\Lambda^\circ} c_\mu \pi (\mu ) = 0$ (as an operator from
$\mathbf{M}^1 (\rd )$ to $\mathbf{M}^\infty (\rd )$), then $c_\mu = 0$
for all $\mu\in\Lambda^\circ$, see~\cite{gr07}. Now, if 
$f = S_{{\boldsymbol\varphi},{\boldsymbol\gamma}} f$ for all $f\in
\mathbf{M}^1(\rd )$, then  by Janssen's representation~\eqref{Eq:MuWijann} we
have  
\[f = s(\Lambda )^{-1} \sum_{\mu\in\Lambda^{\circ}}\sum_{j=1}^n\langle
\varphi_j ,\pi (\mu )\gamma_j\rangle \pi (\mu )f \, .\]
The linear independence of \tfs s implies \eqref{Eq:MuWiWexler}.   The converse is
obvious. \\ 

$(ii)\Leftrightarrow (iii)$: 
  Assume  first that 
$\bigcup_{j=1}^n\mathcal{G} (\varphi_j, \Lambda )$ is a multi-window Gabor frame for $\Lz$. 
The upper bound in \eqref{Eq:RieszCond} follows from the boundedness
of the synthesis operator $\widetilde D_{\g }$ on $\mathbf{L}^2(\rd )$. 
To show the existence  of a lower bound, we apply the Wexler-Raz
relations. Since $\bigcup_{j=1}^n\mathcal{G} (\varphi_j, \Lambda )$ is
a frame with dual  
$\bigcup_{j=1}^n\mathcal{G} (\gamma_j, \Lambda )$  and $\gamma_j\in
\mathbf{M}^1 (\mathbb{R}^d )$ for all $j$,  we have  
$\langle \g , \pi (\mu ) \boldsymbol{\gamma } \rangle =
\sum_{j=1}^n\langle \varphi_j ,\pi (\mu )\gamma_j\rangle =s(\Lambda
)\delta_{\mu ,0}$, and $\mathcal{G} (\g, \Lambda^\circ )$ and
therefore $\mathcal{G} (\boldsymbol\gamma, \Lambda^\circ )$ are
biorthogonal systems in $\mathbf{L}^2 (\mathbb{R}^d,\mathbb{C}^n)$. If
$f = \sum _{\mu \in \Lambda ^\circ  } c_\mu \pi (\mu )\g $, then    
$c_\mu = s(\Lambda )^{-1} \langle \mathbf{f},\pi (\mu
)\boldsymbol\gamma\rangle_{\mathbf{L}^2 (\mathbb{R}^d,\mathbb{C}^n)}$
and  
\[\mathbf{c}= s(\Lambda )^{-1}\widetilde C_{\boldsymbol\varphi,\Lambda^\circ}
\mathbf{f},\]from which the lower bound in \eqref{Eq:RieszCond}
follows.

Conversely, assume  that $\mathcal{G} (\boldsymbol\varphi, \Lambda^\circ )$ is a
Riesz sequence  in $\mathbf{L}^2 (\mathbb{R}^d,\mathbb{C}^n)$. Then
there   exists   a biorthogonal basis of the form  $\{\pi (\mu
)\boldsymbol\gamma :
\mu\in\Lambda^\circ\}$  contained in $\mathcal{K} =
\overline{\spann}(\mathcal{G} (\boldsymbol\varphi, \Lambda^\circ ))$. 
It can be shown that $\boldsymbol{\gamma } \in \mathbf{M}^1(\rd, \bC
^n)$.  The frame
property of $\mathcal{G} (\varphi_j, \Lambda )$ follows from the
Wexler-Raz relations \eqref{Eq:MuWiWexler}.\\ 

With  three classical statements \eqref{Eq:MuWijann} and (ii), (iii)
for  multi-window Gabor frames  the  remaining equivalences follow
exactly as in~\cite{gr07}. 
\end{proof}
  \end{appendix}

\def\cprime{$'$}


\begin{thebibliography}{10}

\bibitem{Berezin71}
F.~A. Berezin.
\newblock Wick and anti-{W}ick symbols of operators.
\newblock {\em Mat. Sb. (N.S.)}, 86(128):578--610, 1971.

\bibitem{BC02}
P.~Boggiatto and E.~Cordero.
\newblock Anti-{W}ick quantization with symbols in {$L\sp p$} spaces.
\newblock {\em Proc. Amer. Math. Soc.}, 130(9):2679--2685 (electronic), 2002.

\bibitem{boch94}
J.-M. {B}ony and J.-Y. {C}hemin.
\newblock {F}unctional spaces associated with the {W}eyl-{H}{\"o}rmander
  calculus. ({E}spaces fonctionnels associ\'es au calcul de
  {W}eyl-{H}{\"o}rmander.).
\newblock {\em {B}ull. {S}oc. {M}ath. {F}rance}, 122(1):77--118, 1994.

\bibitem{BoL}
J.-M. Bony and N.~Lerner.
\newblock Quantification asymptotique et microlocalisations d'ordre
  sup\'erieur. {I}.
\newblock {\em Ann. Sci. \'Ecole Norm. Sup. (4)}, 22(3):377--433, 1989.

\bibitem{brka03}
K.~Brandenburg and M.~Kahrs, editors.
\newblock {\em {A}pplications of Digital Signal Processing to Audio and
  Acoustics}.
\newblock {E}ngineering and Computer science. Kluwer Academic Publishers, 2003.

\bibitem{busby-smith81}
R.~C. Busby and H.~A. Smith.
\newblock Product-convolution operators and mixed-norm spaces.
\newblock {\em Trans. Amer. Math. Soc.}, 263(2):309--341, 1981.

\bibitem{CG03}
E.~Cordero and K.~Gr{\"o}chenig.
\newblock Time-frequency analysis of localization operators.
\newblock {\em J. Funct. Anal.}, 205(1):107--131, 2003.


\bibitem{cogrro06}
E.~{C}ordero, K.~{G}r{\"o}chenig, and L.~{R}odino.
\newblock {L}ocalization operators and time-frequency analysis.
\newblock In N.~e.~a. {C}hong, editor, {\em {H}armonic, {W}avelet and p-adic
  {A}nalysis}, pages 83--109. 2007 {W}orld {S}cient. {P}ubl., 2006.

\bibitem{CF78}
A.~C{\'o}rdoba and C.~Fefferman.
\newblock Wave packets and {F}ourier integral operators.
\newblock {\em Comm. Partial Differential Equations}, 3(11):979--1005, 1978.

\bibitem{Daube88}
I.~Daubechies.
\newblock Time-frequency localization operators: a geometric phase space
  approach.
\newblock {\em IEEE Trans. Inform. Theory}, 34(4):605--612, 1988.


\bibitem{daubechies90}
I.~Daubechies.
\newblock The wavelet transform, time-frequency localization and signal
  analysis.
\newblock {\em IEEE Trans. Inform. Theory}, 36(5):961--1005, 1990.

\bibitem{dofegr06}
M.~{D}{\"o}rfler, H.~G. {F}eichtinger, and K.~{G}r{\"o}chenig.
\newblock {T}ime-{F}requency {P}artitions for the {G}elfand {T}riple $({S}_0,
  \mathbf{L}^2, {S}_0')$.
\newblock {\em {M}ath. {S}cand.}, 98(1):81--96, 2006.

\bibitem{fg89jfa}
H.~G. Feichtinger and K.~Gr{\"o}chenig.
\newblock Banach spaces related to integrable group representations and their
  atomic decompositions. {I}.
\newblock {\em J. Functional Anal.}, 86(2):307--340, 1989.

\bibitem{fg89mh}
H.~G. Feichtinger and K.~Gr{\"o}chenig.
\newblock Banach spaces related to integrable group representations and their
  atomic decompositions. {I}{I}.
\newblock {\em Monatsh. Math.}, 108(2-3):129--148, 1989.

\bibitem{fg92chui}
H.~G. Feichtinger and K.~Gr{\"o}chenig.
\newblock Gabor wavelets and the {H}eisenberg group: {G}abor expansions and
  short time fourier transform from the group theoretical point of view.
\newblock In C.~K. Chui, editor, {\em Wavelets: A tutorial in theory and
  applications}, pages 359--398. Academic Press, Boston, MA, 1992.

\bibitem{fg97jfa}
H.~G. Feichtinger and K.~Gr{\"o}chenig.
\newblock Gabor frames and time-frequency analysis of distributions.
\newblock {\em J. Functional Anal.}, 146(2):464--495, 1997.

\bibitem{FN03}
H.~G. Feichtinger and K.~Nowak.
\newblock A first survey of {G}abor multipliers.
\newblock In {\em Advances in Gabor analysis}, Appl. Numer. Harmon. Anal.,
  pages 99--128. Birkh\"auser Boston, Boston, MA, 2003.

\bibitem{fega06}
C.~{F}ernandez and A.~{G}albis.
\newblock {C}ompactness of time-frequency localization operators on
  ${L}^2({R})$.
\newblock {\em {J}. {F}unct. {A}nal.}, 233(2):335--350, 2006.

\bibitem{FoG05}
M.~Fornasier and K.~Gr{\"o}chenig.
\newblock Intrinsic localization of frames.
\newblock {\em Constr. Approx.}, 22(3):395--415, 2005.

\bibitem{gora98}
S.~J. {G}odsill and P.~J.~W. {R}ayner.
\newblock {\em {D}igital {A}udio {R}estoration.}
\newblock {S}pringer, 1998.

\bibitem{gr01}
K.~{G}r{\"o}chenig.
\newblock {\em {F}oundations of {T}ime-{F}requency {A}nalysis}.
\newblock {A}ppl. {N}umer. {H}armon. {A}nal. {B}irkh{\"a}user {B}oston, 2001.

\bibitem{gr07}
K.~{G}r{\"o}chenig.
\newblock {G}abor frames without inequalities.
\newblock {\em {I}nternational {M}athematics {R}esearch {N}otices 2007},
  2007(article ID rnm111, 21 pages), 2007.

\bibitem{gro07c}
K.~Gr\"ochenig.
\newblock Weight functions in time-frequency analysis.
\newblock In e.~a. L.~Rodino, M.-W.~Wong, editor, {\em Pseudodifferential
  Operators: Partial Differential Equations and Time-Frequency Analysis},
  volume~52, pages 343 -- 366. Fields Institute Comm., 2007.

\bibitem{GL03}
K.~Gr{\"o}chenig and M.~Leinert.
\newblock Wiener's lemma for twisted convolution and {G}abor frames.
\newblock {\em J. Amer. Math. Soc.}, 17:1--18, 2004.

\bibitem{grle04}
K.~{G}r{\"o}chenig and M.~{L}einert.
\newblock {W}iener's {L}emma for {T}wisted {C}onvolution and {G}abor {F}rames.
\newblock {\em {J}. {A}mer. {M}ath. {S}oc.}, 17:1--18, 2004.

\bibitem{GT09}
K.~Gr\"ochenig and J.~Toft.
\newblock Isomorphism properties of {T}oeplitz operators and pseudo-differential
  operators between modulation spaces.
\newblock 2009.
\newblock Preprint.

\bibitem{he03}
C.~{H}eil.
\newblock {A}n introduction to weighted {W}iener amalgams.
\newblock In M.~{K}rishna, R.~{R}adha, and S.~{T}hangavelu, editors, {\em
  {W}avelets and their {A}pplications ({C}hennai, {J}anuary 2002)}, pages
  183--216. {A}llied {P}ublishers, 2003.

\bibitem{jaming98}
P.~Jaming.
\newblock Principe d'incertitude qualitatif et reconstruction de phase pour la
  transform{\'e}e de {W}igner.
\newblock {\em C.~R.~Acad.\ Sci.\ Paris S{\'e}r.\ I Math.}, 327:249--254, 1998.

\bibitem{janssen95}
A.~J. E.~M. Janssen.
\newblock Duality and biorthogonality for {W}eyl-{H}eisenberg frames.
\newblock {\em J. Fourier Anal. Appl.}, 1(4):403--436, 1995.

\bibitem{janssen98}
A.~J. E.~M. Janssen.
\newblock Proof of a conjecture on the supports of {W}igner distributions.
\newblock {\em J. Fourier Anal. Appl.}, 4(6):723--726, 1998.

\bibitem{Lerner}
N.~Lerner.
\newblock The {W}ick calculus of pseudo-differential operators and some of its
  applications.
\newblock {\em Cubo Mat. Educ.}, 5(1):213--236, 2003.


\bibitem{lo07-2}
P.~{L}ouizou.
\newblock {\em {S}peech {E}nhancement: {T}heory and {P}ractice}.
\newblock {C}{R}{C} {P}ress, 2007.



\bibitem{lu09}
F.~{L}uef.
\newblock {P}rojective modules over non-commutative tori are multi-window
  {G}abor frames for modulation spaces.
\newblock {\em {J}. {F}unct. {A}nal.}, 257(6):1921--1946, 2009.

\bibitem{RT93}
J.~Ramanathan and P.~Topiwala.
\newblock Time-frequency localization via the {W}eyl correspondence.
\newblock {\em SIAM J. Math. Anal.}, 24(5):1378--1393, 1993.

\bibitem{ro98}
C.~{R}oads.
\newblock {\em The computer music tutorial}.
\newblock {T}he {M}{I}{T} {P}ress, 1998.

\bibitem{su06}
W.~{S}un.
\newblock {G}-frames and g-{R}iesz bases.
\newblock {\em {J}. {M}ath. {A}nal. {A}ppl.}, 322:437--452, {O}ctober 2006.

\bibitem{Toft04}
J.~Toft.
\newblock Continuity properties for modulation spaces, with applications to
  pseudo-differential calculus. {I}.
\newblock {\em J. Funct. Anal.}, 207(2):399--429, 2004.


\bibitem{wilczock98}
E.~Wilczock.
\newblock Zur {F}unktionalanalysis der {W}avelet- und {G}abortransformation.
\newblock Thesis, TU M{\"u}nchen, 1998.

\bibitem{wo02}
M.-W. {W}ong.
\newblock {\em {W}avelet {T}ransforms and {L}ocalization {O}perators.}
\newblock {O}perator {T}heory: {A}dvances and {A}pplications. 136. {B}asel:
  {B}irkh{\"a}user.  2002.

\end{thebibliography}
\end{document}